\documentclass[12pt]{article}

\usepackage[dvips]{graphicx}
\usepackage{subeqn}
\usepackage[dvipsnames]{color} 
\usepackage{amssymb, amsfonts}

\newtheorem{theorem}{Theorem}
\newtheorem{lemma}{Lemma}
\newtheorem{corollary}{Corollary}

\newtheorem{problem}{Problem}
\newtheorem{remark}{Remark}
\newtheorem{proposition}{Proposition}
\newtheorem{example}{Example}

\newcommand{\eq}{\hspace*{-.2mm}&=&\hspace*{-.2mm}}
\newcommand{\plus}{\hspace*{-.2mm}&+&\hspace*{-.2mm}}

\newcommand{\eps}{\varepsilon}
\newcommand{\R}{{\mathbb R}}

\newcommand{\kb}{{\bf k}}

\newcommand{\Dc}{{\cal D}}

\def\calf{\mathcal F}

\title{On surfaces with a prescribed curvilinear projection \\ of one field of principal directions}

\author{Vladimir Rovenski and Leonid Zelenko \\
         \small Department of Mathematics,
             Faculty of Science and Science Education\\
          \small University of Haifa, Mount Carmel, Haifa, 31905, Israel\\          
          \small E-mail: rovenski@math.haifa.ac.il
                                           and zelenko@math.haifa.ac.il}

\date{}

\begin{document}

\maketitle

\begin{abstract}
A class of surfaces-graphs in a Riemannian 3-space with a prescribed projection of one field of principal directions onto a surface $\Pi$ is considered.
A~problem of determination of such surfaces when both principal curvatures are given over a line in~$\Pi$ is formulated and studied.
The geometric problem is reduced to the Cauchy problem for quasilinear
PDE's which, under certain conditions for data, are hyperbolic and admit a unique solution.
It is shown that the parallel curved (PC) surfaces in space forms
provide a special class of global solutions to the geometrical problem
with weaker regularity assumptions. Such solutions may be found by an iteration function sequence.


\end{abstract}


\section*{Introduction}
\addcontentsline{toc}{section}{Introduction}
\label{S-01}

 The surfaces possessing nontrivial deformations which preserve
principal curvatures and directions (or, equivalently, the shape operator) were
investigated by several authors,
see \cite{Bryant_2001}, \cite{V} and review with bibliography in~\cite{Fer}.
It is known that surfaces with one family of principal curves being geodesic
(as for parallel curved (PC) surfaces recently studied in \cite{Ando2004}\,--\,\cite{Ando2007})
represent degenerate case in studying immersions of simply connected surfaces with a prescribed shape operator.
Recent studying of reconstruction of surfaces by their partially given  principal curvatures and directions may be useful for applications of differential geometry to computer graphics, the wavefront analysis in applied optics, etc.

In what follows, $(\bar M^3,\bar g)$ denotes a $C^3$-regular Riemannian 3-space with
coordinates $x_1,x_2$, $x_3$ ($|x_i|\le a_i$) for some $a_i\in\R$,
$\Pi=\{x_3=0\}$ a $C^3$-regular surface,
$\gamma=\{x_2=x_3=0\}$ the coordinate curve,
and $\pi(x_1,x_2,x_3)=(x_1,x_2,0)$ the \textit{curvilinear projection}.

In the paper we consider surfaces-graphs $M^2\subset\bar M^3$ with \textit{prescribed curvilinear projection (onto $\Pi$) of one field of principal directions}.
We show (see Theorem~\ref{T-main0}) that such surfaces depend on two arbitrary functions of one variable,
namely, the principal curvatures over $\gamma$ which are assumed close enough to corresponding values for~$\Pi$.
More precisely, we study the following.

\begin{problem}\label{Prob-1}\rm
Given $(\bar M^3,\bar g)$, a vector field $l$ transversal to $\gamma$ on $\Pi$,
functions $\bar k_1$ and $\bar k_2$ of class $C^0(\gamma)$,
\textbf{find} a function $f$ of class $C^2(\Pi)$, whose
graph $M^2: x_3=f(x_1,x_2)$ in $\bar M^3$ satisfies the conditions:

(i) \textit{the projection ($\pi$) onto $\Pi$ of the field $\partial_1$ of principal directions}

 \textit{corresponding to the principal curvature $k_1$ (of $M^2$) coincides with $l$},

(ii) \textit{the principal curvatures $k_i$ (of $M^2$) over $\gamma$
 coincide with $\bar k_i$: $k_{i|\,\gamma}=\bar k_i$},

(iii) \textit{the values of $f$ and $df$ at the point $(0,0,0)$ of $\gamma$ are given}.
\end{problem}

Our approach is based on reducing the \textbf{Problem~\ref{Prob-1}} to the Cauchy problem
for a quasilinear system of PDE's which, under certain conditions for data, is hyperbolic and admits a unique local smooth solution.
 The PC surfaces in $\R^3$
 represent a special class of solutions when a family of curvature lines projects onto $\Pi$ as parallel lines or concentric circles. Such surfaces are recovered by an iteration function sequence and using the reconstruction of \underline{two}
planar curves by their curvature (see Theorem~\ref{T-RPC2} and Proposition~\ref{T-RPC1}).
Notice that the space of PC surfaces free of umbilics and having the same shape operator depends on \underline{one} arbitrary function of one variable, see \cite{Bryant_2001}.

The structure of the work is the following.
\textbf{Section~\ref{sec:main-res}} represents main results
(Theorems~\ref{T-main0}\,--\,\ref{T-RPC2} and Corollary~\ref{T-main3} for $M^2\subset\R^3$).
\textbf{Section~\ref{sec:proofs}} contains proofs.
\textbf{Section~\ref{sec:appendix}} contains necessary facts on PC surfaces and examples.

\section{Main results}
\label{sec:main-res}

We shall use the following notation:
$u=(u_1,\dots,u_n)\in\R^n$; $\|u\|_\infty=\max_{1\le i\le n}|u_i|$,
$C^0(\Dc)$ the linear space of bounded continuous functions
$u:\Dc\rightarrow\R^m$ ($\Dc$ is a domain in $\R^n$);
$\|u\|_{\Dc}=\sup_{x\in\Dc}\|u(x)\|_\infty$ the norm in $C^0(\Dc)$;
$C^k(\Dc)$ is the set of functions
$u:\Dc\rightarrow\R^m$, having in $\Dc$ continuous partial derivatives of order~$k$.
For short, we omit $m$ from the above notations.

For simplicity, we assume in what follows that
$\Pi$ is a totally umbilical surface with the normal curvature $\lambda$
(if $\lambda\equiv0$ then $\Pi$ is totally geodesic).

 The main result of the paper is the following.

\begin{theorem}\label{T-main0}
Let $\bar k_{1},\bar k_{2}$ be functions of class $C^{1}(\gamma)$,
$l$ a vector field of class $C^{2}(\Pi)$,
that is transversal but not orthogonal to~$\gamma$.
If $\|\bar k_{i}-\lambda_{\,|\,\gamma}\|_{\,\gamma}$ are small
enough,
then \textbf{Problem~\ref{Prob-1}}
admits in $(\bar M^3,\bar g)$ a unique
local solution.
 Namely, there are $\Delta,K>0$ such that if
 $\|\bar k_{i}{-}\lambda_{\,|\,\gamma}\|_\gamma{<}\,\Delta$ $(i=1,2)$,
 then for some $\eps\in(0,a_2]$ there \underline{exists} a
 function $f$ of class $C^{3}$ on
 $\Pi_{K,\eps}=\{|x_1|+K x_2\le a_1,\,0\le x_2\le\eps,\,x_3=0\}$
 with the properties:
the~principal curvatures $k_i$ of $M^2:\,x_3=f(x_1,x_2)$ satisfy
$k_{i|\,\gamma}=\bar k_{i}$,
$l$ is tangent to the $\pi$-projection onto $\Pi_{K,\eps}$ of $k_1$-curvature lines,
and $f(0,0)=df(0,0)=0$.
 Moreover, there is $r\in(0,a_3]$ such that the solution $f$
 is \underline{unique} in the class of $C^3$-regular
 functions satisfying $\|(f, f_{x_1}, f_{x_2})\|_{\Pi_{K,\eps}}\le r$.
\end{theorem}

 One may apply Theorem~\ref{T-main0} to surfaces in 3-space forms (see also
 Section~\ref{sec:warped2}). We illustrate this for $\R^3$ with cartesian coordinates.

\begin{corollary}\label{T-main3}
Let $\alpha\ne0$ be a function of class $C^{2}$ on a rectangle $\Pi=\{|x|\le a_1,\,|y|\le a_2,\,z=0\}$ in $\R^3$ (with cartesian coordinates), and $\bar k_{1},\bar k_{2}$ functions
of class $C^{1}$ on the segment $\gamma=\{|x|\le a_1,\,y=z=0\}$.
If $\|\bar k_{i}\|_{\,\gamma}$ are small enough, then \textbf{Problem~\ref{Prob-1}} admits
a unique
local solution. Namely, there are $\Delta,K>0$ such that
if $\|\bar k_{i}\|_\gamma<\Delta$, then for some $\eps\in(0,a_2]$ there \underline{exists}
a function $f$ of class $C^{3}$ on $\Pi_{K,\eps}=\{|x|+K y\le a_1,\,0\le y\le\eps,\,z=0\}$
with the properties:
the principal curvatures $k_i$ of the graph $M^2: z=f(x,y)$ in $\R^3$ satisfy
$k_i(x,0)= \bar k_{i}(x)$, the vector field $\alpha\,\partial_x+\partial_y$ is tangent
to the projection onto $\Pi_{K,\eps}$ of $k_1$-curvature~lines,
and $f(0,0)=df(0,0)=0$.
 Moreover, there is $r>0$ such that the solution $f$ is \underline{unique} in the class of $C^3$-regular functions satisfying $\|(f, f_{x}, f_{y})\|_{\Pi_{K,\eps}}\le r$.
\end{corollary}

\begin{remark}\label{R-01}\rm
 (a)
 The condition that $\Pi$ is totally umbilical can be dropped.
In this case, $\bar k_i$ should be close enough to corresponding principal curvatures of $\Pi$ along $\gamma$, and $l$ close enough to one of principal directions on $\Pi$.
The values of $f(0,0)$ and $df(0,0)$ can be taken small enough (see Proposition~\ref{P-twoODEs}), in the present text for simplicity we assume them zero.

 If $\alpha=\mbox{const}\ne0$ in Corollary~\ref{T-main3}, then $M^2$ is a PC surface,
 i.e., the planes $\{x-\alpha\,y=c\}$ intersect $M^2$ by curvature lines, see Section~\ref{sec:Ex-003}. PC surfaces in spherical coordinates, see also Section~\ref{sec:Ex-003},  illustrate Theorem~\ref{T-main0}.

 (b) A \textit{normal geodesic graph in $(\bar M^3,\bar g)$ of a function} $f:\Pi\to\R,\ |f|<r_{f}(\Pi)$ ($r_{f}(\Pi)$ is the focal radius of $\Pi$)
is a surface $M^2=\bigcup_{x\in\Pi}\{\gamma_x(f(x))\}\subset\bar M^3$, where
$\gamma_x(t)\ (x\in\Pi)$ is a unit speed geodesic normal to $\Pi$.

The \textit{semi-geodesic coordinates} $(x_1,x_2)$ on $\Pi$ with the base curve $\{x_2=0\}$; have the metric is $g_2=\xi^2(x_1,x_2)\,dx_1^2+dx_2^2$, where
 $\xi_{,22}+K(x_1,x_2)\xi=0$, $\xi(x_1,0)=1,\ \xi_{,2}(x_1,0)=0$ (see \cite{Top}),
and $K$ is the gaussian curvature of $\Pi$.

Let $\bar M^3$ has the coordinates $x_1,x_2,x_3$ ($|x_i|\le a_i$) such that

\textbf{--} $(x_1,x_2)$ are semi-geodesic coordinates on a surface $\Pi=\{x_3=0\}$,

\textbf{--} the curve $\gamma=\{x_2=x_3=0\}$ is a simple geodesic in $\Pi$, and

\textbf{--} $x_3$ is the signed distance to $\Pi$
(hence $\bar g_{i3}=\delta_{i3}$ and $\bar\Gamma^k_{33}=0$).

\noindent
 One may obtain corollary of Theorem~\ref{T-main0},
 where $\pi:\bar M^3\to\Pi$ means ``the nearest point" in~$\Pi$
 (for $\R^3$ we again have Corollary~\ref{T-main3}).
\end{remark}

We will formulate for $\bar M^3$ (and study in the paper for $\R^3$) the problem that plays essential role in solving \textbf{Problem~\ref{Prob-1}} for PC surfaces.

\begin{problem}\label{Prob-2}\rm
Given
 a vector field $l=l(x_1,x_2)$ on $\Pi\subset\bar M^3$,
 an integral curve $\gamma_1\subset\Pi$ of $l$ through $O\in\Pi$,
 a curve $\gamma_2\subset\Pi$ transversal to $\gamma_1$ through $O$
and
 functions $\bar k_i\in C^{0}(\gamma_i)$ on $\gamma_i$ ($i=1,2$),
\textbf{find} a function $f$ of class $C^2(\Pi)$, whose
graph $M^2: x_3=f(x_1,x_2)$ in $\bar M^3$ satisfies the conditions:

(i) \textit{the projection ($\pi$) onto $\Pi$ of the field $\partial_1$ of principal directions}

 \textit{corresponding to the principal curvature $k_1$ of $M^2$, coincides with $l$},

(ii) \textit{the principal curvatures $k_{i}$ of $M^2$ coincide with
$\bar k_i$ over $\gamma_i$: $k_{i|\,\gamma_i}=\bar k_i$},

(iii) \textit{the values of $f$ and $df$ at the point $O$ are given}.
\end{problem}

Define a rectangle $\Pi(a)=\{(x,y): |\alpha x+y|\le a\,\alpha,\,|x-\alpha y|\le a\}$ in the $xy$-plane of~$\R^3$, where $a,\alpha\in\R$ are positive.
The PC surfaces represent a special class of global solutions (i.e., on the domains $\Pi(a)$ with an arbitrary  $a<a_1$) to
\textbf{Problems~\ref{Prob-1}} and \textbf{\ref{Prob-2}} with weaker regularity assumptions.

\begin{theorem}\label{T-RPC2}
\hskip-0.8pt Let $a_1,\alpha>0$ be real, $\bar k_{1}\in C^{1}([-a_1,a_1])$
and $\bar k_{2}\in C^{0}([-a_1,a_1])$.
Then for any $a\in(0,a_1)$ one can choose $\delta>0$  such that
if $\|\bar k_i\|_{[-a_1,a_1]}<\delta$ $(i=1,2)$, then there
exists a unique PC surface $M^2\subset\R^3: z=f(x,y)$, where
$f\in C^2(\Pi(a))$, with the properties

\hskip-2mm
 -- the $k_1$-principal direction $\partial_1$ on $M^2$
 is parallel to the plane $x-\alpha\,y=0$,

\hskip-2mm
 -- the principal curvatures of $M^2$ satisfy $k_i(x,0)=\bar k_{i}(x)\ (|x|\le a,\,i=1,2)$.
\end{theorem}

Theorem~\ref{T-RPC2} is based on the following result concerning existence of a solution to \textbf{Problem~\ref{Prob-2}} in the class of PC surfaces.

\begin{proposition}\label{T-RPC1}
 Given $a_1,\alpha>0$, let
 $\gamma_1=\{(\alpha\,u, u):\,|u|\le\frac{a_1\alpha}{\sqrt{\alpha^2+1}}\}$ and
 $\gamma_2=\{(x, 0):\,|x|\le a_1\}$ be line segments in $xy$-plane of $\R^3$.
 Let $\tilde k\in C^1(\gamma_1)$ and $\bar k_2\in C^0(\gamma_2)$.
Then for any $a\in(0,a_1)$ there is $\delta>0$ such that~if
$
  \bar k=\max\limits_{|x|\,\le\,a_1}\{\,
  |\tilde k(x\,\alpha/\sqrt{\alpha^2+1})|,\ |\bar k_2(x)|\,\}<\delta,
$
 then there exists a unique PC surface $M^2: z=f(x,y)$ in $\R^3$, where
$f\in C^2(\Pi(a))$, with the properties

\hskip-3mm
 -- the principal direction $\partial_1$ is orthogonal to the vector $e_1-\alpha e_2$,

\hskip-3mm
 -- $\tilde k$ is the principal curvature of $M^2$ over
 $\gamma_1$ corresponding to $\partial_1$, and

\hskip-3mm
 -- $\bar k_2$ is the principal curvature of $M^2$ over $\gamma_2\cap\Pi(a)$ corresponding to

 the second principal direction $\partial_2$.

\noindent
The solution can be found using reconstruction of two planar curves by their curvature functions.
\end{proposition}

\section{Proofs}
\label{sec:proofs}

In~{Section~\ref{sec:results}} we prove Theorem~\ref{T-main0} and its Corollary~\ref{T-main3}.
Based on the Euler formula for the principal curvatures we deduce a system of PDE's (Proposition~\ref{T-compR}).
Using compatibility conditions, we transform above equations to equivalent quasi-linear system (Proposition~\ref{T-01}), for which we formulate the Cauchy Problem.
The initial values are analyzed in Lemma~\ref{L-04}, where $M^2$ is recovered over $\gamma$.
{Section~\ref{sec:Th2-3}} shows that the PC surfaces in space forms
(i.e., a family of curvature lines projects onto parallel lines or concentric circles)
represent a special class of global solutions to the geometrical problem with weaker regularity assumptions
($\bar k_{1}\in C^{1}([-a_1, a_1]),\bar k_{2}\in C^{0}([-a_1, a_1])$).
We~approximate such solutions by the iterated function
sequence, and apply the procedure of reconstruction of two plane curves by their curvature functions (Theorem~\ref{T-RPC2} and Proposition~\ref{T-RPC2-B}).
{Section~\ref{sec:axresults}} contains auxiliary lemmas
and necessary facts on hyperbolic PDE's.

\subsection{Proof of Theorem~\ref{T-main0} and its corollary}
 \label{sec:results}

A \textit{surface-graph} $M^2$ in $\bar M^3$ is defined by equation $x_3=f(x_1,x_2)$.
Let $\bar e_i=\frac{\partial}{\partial x_i}$ be coordinate vector fields on $\bar M^3$,
$\hat e_i$ their restrictions on $M^2$, and $e_1,e_2$ the coordinate vector fields
(lifts under $\pi$ of $\bar e_1,\bar e_2$ from $\Pi$) on~$M^2$.
For~simplicity assume in what follows that $\bar g_{13}=\bar g_{23}=0$
(see also Remark~\ref{R-01}).
We have $e_1=\hat e_1+p\,\hat e_3$ and $e_2=\hat e_2+q\,\hat e_3$,
where $p=f_{x_1}$ and $q=f_{x_2}$.
The~metric on $M^2$ is given by $g_{ij}=\bar g(e_i,e_j)$.
The~coefficients of the first and the second fundamental forms of $M^2$ are denoted by
\begin{eqnarray*}
 E=g_{11}(x_1,x_2),\quad
 F=g_{12}(x_1,x_2)=g_{21}(x_1,x_2),\quad
 G=g_{22}(x_1,x_2),\\
 L= b_{11}(x_1,x_2),\quad
 M=b_{12}(x_1,x_2)=b_{21}(x_1,x_2),\quad
 N=b_{22}(x_1,x_2).
\end{eqnarray*}
Suppose that $\partial_1=\alpha\,e_1+\,e_2,\ \partial_2=\beta_1\,e_1+\beta_2\,e_2$
are the principal directions on $M^2$. We compute the product
$0=\bar g(\partial_1,\partial_2)=\beta_1(\alpha E+F)+\beta_2(\alpha F+G)$.
Hence $\beta_1:\beta_2=-(\alpha F+G):(\alpha E+F)$.
Denote by $\hat g_{ij}=\bar g(\hat e_i,\hat e_j)$.

From above and Lemma~\ref{L-EFG-1} (Section~\ref{sec:axresults}) it follows
\begin{eqnarray}\label{E-alpha-3}
\nonumber
 \hskip-6mm&&
 \alpha^2 E+2\,\alpha F+G = g(\partial_1,\partial_1)>0,\\
 \hskip-6mm&&\alpha E+F = \hat g_{33}(\alpha p+q)p +\alpha\hat g_{11}+\hat g_{12},\\
\nonumber
 \hskip-6mm&& EG-F^2 = \hat g_{33}(\hat g_{22}p^2+\hat g_{11}q^2-2\hat g_{12}p\,q)
 +\hat g_{11}\hat g_{22}-\hat g_{12}^2\ge\hat g_{11}\hat g_{22}-\hat g_{12}^2.
\end{eqnarray}
Define the functions $H_{ij}=H_{ij}^{(1)}k_1+H_{ij}^{(2)}k_2+H_{ij}^{(0)}\ (i,j=1,2)$, where
\begin{equation}\label{E-HijkR}
\hskip-4mm
\begin{array}{c}
 H_{11}^{(1)}=\frac{\delta(\alpha E+F)^2}{\alpha^2 E+2\,\alpha F+G},\ \
 H_{11}^{(2)}=\frac{\delta(EG-F^2)}{\alpha^2 E+2\,\alpha F+G},\
 H_{12}^{(1)}=H_{21}^{(1)}=\delta\frac{(\alpha E+F)(\alpha F+G)}{\alpha^2 E+2\,\alpha F+G},\\
 H_{12}^{(2)}=H_{21}^{(2)}=-\frac{\delta\,\alpha(EG-F^2)}{\alpha^2 E+2\,\alpha F+G},\
 H_{22}^{(1)}=\frac{\delta(\alpha F+G)^2}{\alpha^2 E+2\,\alpha F+G},\
 H_{22}^{(2)}=\frac{\delta\,\alpha^2(EG-F^2)}{\alpha^2 E+2\,\alpha F+G},\\
 H_{11}^{(0)}=-L_1,\quad
 H_{12}^{(0)}=-M_1,\quad
 H_{22}^{(0)}=-N_1,
\end{array}
\end{equation}
$\delta=\sqrt{(EG-F^2)/\det\hat g}\ge 1/\sqrt{\hat g_{33}}$, see (\ref{E-alpha-3})$_3$,
and $L_1,M_1,N_1$ are given in (\ref{E-L1-M1-N1}) of Section~\ref{sec:axresults}.
For $\R^3$ we get $H_{ij}=H_{ij}^{(1)}k_1+H_{ij}^{(2)}k_2$, where
$H_{21}^{(i)}=H_{12}^{(i)}$,
\begin{equation}\label{E-HijkR-R3}
\hskip-12mm
\begin{array}{ccc}
 && H_{11}^{(1)}= \frac1{\delta_1}{\sqrt{1{+}{p}^{2}{+}{q}^{2}}
   ({\alpha}({p}^{2}{+}1)+{p}\,{q})^2},\quad
   H_{11}^{(2)}=\frac1{\delta_1}(1{+}{p}^{2}{+}{q}^{2})^{\frac32},\\
 &&
 H_{12}^{(1)}{=}\frac1{\delta_1}{\sqrt{1{+}{p}^{2}{+}{q}^{2}}
   ({\alpha}p\,q{+}q^2{+}1)(\alpha(p^2{+}1){+}p\,q)},\
   H_{12}^{(2)}{=}-\frac\alpha{\delta_1}(1{+}{p}^{2}{+}{q}^{2})^{\frac32},\\
 && H_{22}^{(1)}=\frac1{\delta_1}\sqrt{1{+}{p}^{2}{+}{q}^{2}}
   ({\alpha}{p}\,{q}{+}{q}^{2}{+}1)^2,\quad
   H_{22}^{(2)}=\frac{\alpha^{2}}{\delta_1}(1{+}{p}^{2}{+}{q}^{2})^{\frac32},
\end{array}
\end{equation}
and $\delta_1=(\alpha p+q)^2+\alpha^2+1\ge 1$.

\begin{example}\label{Ex-Pi}\rm
Let $M^2=\Pi=\{x_3=0\}$, hence $f=p=q=0$.
By Lemma~\ref{L-EFG-1} (Section~\ref{sec:axresults}),
the coefficients of the 1-st and the 2-nd fundamental forms are
\[
 \hat E=\bar g_{11},\ \hat F=\bar g_{12},\ \hat G=\bar g_{22},\quad
 \hat L=\bar\Gamma^3_{11}\sqrt{\bar g_{33}},\
 \hat M=\bar\Gamma^3_{12}\sqrt{\bar g_{33}},\
 \hat N=\bar\Gamma^3_{22}\sqrt{\bar g_{33}}.
\]
The functions (\ref{E-HijkR}) on $\Pi$ have the following form:
\[
\begin{array}{ccc}
 &&\hskip-9mm \hat H_{11}^{(1)}=\delta\frac{(\alpha\bar g_{11}+\bar g_{12})^2}
 {\alpha^2\bar g_{11}{+}2\,\alpha \bar g_{12}{+}\bar g_{22}},\
 \hat H_{11}^{(2)}=\delta\frac{\bar g_{11}\bar g_{22}-\bar g_{12}^2}
 {\alpha^2\bar g_{11}{+}2\,\alpha\bar g_{12}{+}\bar g_{22}},\
 \hat H_{12}^{(1)}=
 \delta\frac{(\alpha\bar g_{11}{+}\bar g_{12})(\alpha\bar g_{12}{+}\bar g_{22})}
 {\alpha^2\bar g_{11}+2\,\alpha\bar g_{12}+\bar g_{22}},\\
 &&\hskip-8mm \hat H_{12}^{(2)}=
 -\delta\frac{\alpha(\bar g_{11}\bar g_{22}-\bar g_{12}^2)}
 {\alpha^2\bar g_{11}+2\,\alpha\bar g_{12}+\bar g_{22}},\
 \hat H_{22}^{(1)}=\delta\frac{(\alpha\bar g_{12}+\bar g_{22})^2}{\alpha^2\bar g_{11}+2\,\alpha \bar g_{12}+\bar g_{22}},\
 \hat H_{22}^{(2)}=\delta\frac{\alpha^2(\bar g_{11}\bar g_{22}-\bar g_{12}^2)}{\alpha^2\bar g_{11}+2\,\alpha \bar g_{12}+\bar g_{22}},\\
 &&\hskip-8mm \hat H_{11}^{(0)}=-\bar\Gamma^3_{11}\sqrt{\bar g_{33}},\qquad
 \hat H_{12}^{(0)}=-\bar\Gamma^3_{12}\sqrt{\bar g_{33}},\qquad
 \hat H_{22}^{(0)}=-\bar\Gamma^3_{22}\sqrt{\bar g_{33}},
\end{array}
\]
where $\delta=1/\sqrt{\bar g_{33}}$.
The function $\hat\alpha$ (of the principal direction of $\Pi$) satisfies
\begin{equation}\label{E-tumb1}
 (\bar g_{11}\bar\Gamma^3_{12}-\bar g_{12}\bar\Gamma^3_{11})\,\hat\alpha^2
 +(\bar g_{11}\bar\Gamma^3_{22}-\bar g_{22}\bar\Gamma^3_{11})\,\hat\alpha
 +(\bar g_{12}\bar\Gamma^3_{22}-\bar g_{22}\bar\Gamma^3_{12})=0.
\end{equation}
The principal curvatures $\hat k_i$ of $\Pi$ are solutions to the quadratic equation
\begin{equation}\label{E-tumb2}
 \hskip-2mm(\bar g_{11}\bar g_{22}{-}\bar g_{12}^2)\/k^2
 {-}\sqrt{\bar g_{33}}\/(\bar g_{11}\bar\Gamma^3_{22}{+}\bar g_{22}\bar\Gamma^3_{11}
 {-}2\bar g_{12}\bar\Gamma^3_{12})\/k
 {+}\bar g_{33}(\bar\Gamma^3_{11}\bar\Gamma^3_{22}{-}(\bar\Gamma^3_{12})^2)\,{=}\,0.
\end{equation}
Since $\Pi$ is totally umbilical, $\hat k_1=\hat k_2=\lambda$ are the roots of (\ref{E-tumb2}), hence
\begin{equation}\label{E-tumb4}
 \bar\Gamma^3_{ij}=-\lambda\bar g_{ij}/\bar g_{33}\quad(1\le i,j\le3).
\end{equation}
In this case, (\ref{E-tumb1}) is satisfied by any $\hat\alpha$.
\end{example}

First, we will prove Propositions~\ref{T-compR}, \ref{T-01} and Lemma~\ref{L-04}.

\begin{proposition}\label{T-compR}
Let $M^2\subset\bar M^3$ be the graph of $f\in C^2(\Pi)$.
Then $k_1,k_2\in C^{0}(\Pi)$ are the principal curvatures
and $l=\alpha\,\bar e_1+\bar e_2$ (with $\alpha\in C^1(\Pi)$) is the projection onto $\Pi$ of $k_1$-principal direction of $M^2$ if and only~if
\begin{subequations}
\begin{eqnarray}
\label{E-02-R}
 &&\hskip-14mm p_{x_1} = H_{11}(x_1,x_2, f,p,q,k_1,k_2),\
 q_{x_1}{=} H_{12}(x_1,x_2, f,p,q, k_1,k_2),\ f_{x_1}{=}\,p\\
\label{E-04-R}
 &&\hskip-14mm p_{x_2} = H_{21}(x_1,x_2, f,p,q,k_1,k_2),\
 q_{x_2}{=} H_{22}(x_1,x_2, f,p,q,k_1,k_2),\ f_{x_2}{=}\,q.
\end{eqnarray}
\end{subequations}
If $f\in C^3(\Pi)$, $k_1,k_2\in C^{1}(\Pi)$ and $\alpha E+F\ne0$, then
the compatibility conditions for (\ref{E-02-R},b) are reduced to PDE's
\begin{equation}\label{E-3.2-R}
\begin{array}{cc}
 k_{1,x_2}-\frac{\alpha F{+}G}{\alpha E+F}\,k_{1,x_1}=\Psi_1,\qquad
 k_{2,x_2}+\alpha\,k_{2,x_1} =\Psi_2,
\end{array}
\end{equation}
where $\Psi_i(x_1,x_2, f,p,q,k_1,k_2)$ are known functions (see the proof).
The characteristics of (\ref{E-3.2-R}) are the projections onto $\Pi$
of curvature lines of $M^2$.
\end{proposition}

\hskip-3mm
\textbf{Proof}. Let $V{=}\mu_1\,e_1{+}\mu_2\,e_2$ be a vector on $M^2$.
The functions $k_1,k_2{\in} C^{0}(\Pi)$ are the principal curvatures of $M^2$,
and $l=\alpha\,e_1+e_2$ is the projection onto $\Pi$ of $k_1$-principal direction
if and only if they satisfy the Euler formula
 $k_n(V)=[\frac{g(V,\partial_1)^2}{I(\partial_1)}\,k_1
 {+}\frac{g(V,\partial_2)^2}{I(\partial_2)}\,k_2]/I(V)$,
with the normal curvature $k_n(V)=\frac{II(V)}{I(V)}$. Hence
\begin{equation}\label{E-euler}
 L\mu_1^2+2M\mu_1\mu_2+N\mu_2^2=
 [\,g(V,\partial_1)^2/I(\partial_1)]k_1+[\,g(V,\partial_2)^2/I(\partial_2)]k_2.
\end{equation}
Notice that if $\partial_1=\alpha\,e_1+e_2$ is the principal direction on $M^2$
then $\partial_2=-(\alpha F+G) e_1+(\alpha E+F) e_2$ is the second principal direction. We have
\begin{eqnarray*}
 I(\partial_1)\eq\alpha^2 E+2\,\alpha F+G,\qquad
 I(\partial_2)=(EG-F^2)(\alpha^2 E+2\,\alpha F+G),\\
 g(V,\partial_1)\eq \alpha \mu_1E+(\alpha\mu_2+\mu_1)F+\mu_2G
                  =(\alpha E+F)\mu_1 + (\alpha F+G)\mu_2,\\
 g(V,\partial_2)\eq (EG-F^2)(-\mu_1+\alpha\mu_2),\quad
 I(V)=E\mu_1^2+2F\mu_1\mu_2+G\mu_2^2.
\end{eqnarray*}
Since $V$ is arbitrary, (\ref{E-euler}) is equivalent to the system
\[
 \delta\,L = H_{11}^{(1)}\,k_1 + H_{11}^{(2)}\,k_2,\ \
 \delta\,M = H_{12}^{(1)}\,k_1 + H_{12}^{(2)}\,k_2,\ \
 \delta\,N = H_{22}^{(1)}\,k_1 + H_{22}^{(2)}\,k_2,
\]
which, in view of
 $p_{x_1}=\delta\,L-L_1,\,p_{x_2}=\delta\,N-N_1$ and $p_{x_2}=q_{x_1}=\delta\,M-M_1$
 (see Lemma~\ref{L-EFG-1} in Section~\ref{sec:axresults}), yields (\ref{E-02-R},b).

If $f\in C^3(\Pi)$ and $k_1,k_2\in C^{1}(\Pi)$,
then the compatibility conditions for (\ref{E-02-R},b), i.e.,
$
 (p_{x_1})_{x_2}=(p_{x_2})_{x_1},\ (q_{x_1})_{x_2}=(q_{x_2})_{x_1},
$
take a form
\begin{equation}\label{E-18-R}
\begin{array}{c}
 H_{11,x_2} {+} H_{11,\,p} H_{12}{+}H_{11,\,q} H_{22}{+}H_{11,f} q =
 H_{12,x_1} {+} H_{12,\,p} H_{11}{+}H_{12,\,q} H_{12}{+}H_{12,f} p,\\
 H_{12,x_2} {+} H_{12,\,p}H_{12}{+}H_{12,\,q}H_{22}{+}H_{12,f} q =
 H_{22,x_1} {+} H_{22,\,p}H_{11}{+}H_{22,\,q}H_{12}{+}H_{22,f} p.
\end{array}
\end{equation}
Substituting (\ref{E-02-R},b) into (\ref{E-18-R}), we obtain PDE's
for $\kb=(k_1, k_2)$
\begin{equation}\label{E-3.1-R}
 A(x_1,x_2, f,p,q)\,\kb_{,x_2} + B(x_1,x_2, f,p,q)\,\kb_{,x_1}+b(x_1,x_2, f,p,q, \kb)=0,
\end{equation}
where the matrices
$A=\Big(\begin{array}{cc}
 -H_{11}^{(1)}& -H_{11}^{(2)}\\
  H_{12}^{(1)}& H_{12}^{(2)}
 \end{array}\Big)$,
 $B=\Big(\begin{array}{cc}
  H_{12}^{(1)}& H_{12}^{(2)}\\
 -H_{22}^{(1)}&-H_{22}^{(2)}
\end{array}\Big)$
are $C^1$-regular, and the components of the vector $b=(b_1,b_2)$ are
\begin{eqnarray*}
 b_1\eq(H^{(1)}_{12,x_1}{-}H^{(1)}_{11,x_2})k_1{+}(H^{(2)}_{12,x_1}{-}H^{(2)}_{11,x_2})k_2\\
 &&+H_{12,\,p} H_{11}{+}H_{12,\,q} H_{12}-H_{11,\,p} H_{12}{-}H_{11,\,q} H_{22}
 +H_{12,f} p-H_{11,f} q,\\
 b_2\eq(H^{(1)}_{12,x_2}{-}H^{(1)}_{22,x_1})k_1{+}(H^{(2)}_{12,x_2}{-}H^{(2)}_{22,x_1})k_2\\
 &&+H_{12,\,p}H_{12}{+}H_{12,\,q}H_{22}-H_{22,\,p}H_{11}{-}H_{22,\,q}H_{12}
 +H_{22,f} p -H_{12,f} q.
\end{eqnarray*}
 Notice that $\det A=\delta^2\frac{(EG-F^2)(\alpha E+F)}{\alpha^2 E+2\alpha F+G}\ne0$ when
 $\alpha E+F\neq 0$.
 The direct computation shows that
$
 A^{-1}B=\Big(\begin{array}{cc}
 \lambda_1 & 0\\ 0 & \lambda_2 \\
 \end{array}\Big)
$,
where  $\lambda_1{=}-\frac{\alpha F+G}{\alpha E+F}$ and $\lambda_2=\alpha$.
Hence (\ref{E-3.1-R}) is equivalent to the system (\ref{E-3.2-R})
with $(\Psi_1,\Psi_2)^T=A^{-1} b$.\hfill$\square$

\begin{remark}\label{R-petcod}\rm
Let $M^2\subset{\mathbb R}^3$ be a $C^3$-regular surface without umbilical points parameterized by coordinates $u,v$ of the curvature lines, ${\bf r}\,(u,v)$. Hence $F=M=0$.
Let $k_1,k_2$ be the principal curvatures of $M^2$.
 In this case, if $\alpha=0$, then $H^{(1)}_{22}=\delta\,G$,
 $H^{(2)}_{11}=\delta E$ and other $H^{(k)}_{ij}$ are zero.
 Hence $b_{21}=\delta\,G$, $a_{12}=\delta E$ and other $a_{ij}, b_{ij}$ are zero.
 Hence, (\ref{E-3.1-R}) is reduced to the system
 $\delta E k_{1,v}+\tilde c_1=0,\ \delta G k_{2,u}+\tilde c_2=0$
 that is equivalent to Peterson-Codazzi formulae (see \cite{Top})
 $k_{1,v}=(k_2-k_1)\frac{E_{,v}}{2 E},\
 k_{2,u}=(k_1-k_2)\frac{G_{,u}}{2 G}$.
 One may compare the last formulae with our compatibility equations~(\ref{E-18-R}).
\end{remark}

\begin{proposition}\label{T-01}
Let the functions $f, p, q, k_1,k_2,\alpha$ of class $C^1$ on $\Pi$ (in particular, on $\Pi_{K,\eps}$) satisfy (\ref{E-04-R}).

(i) If (\ref{E-02-R}) holds in~$\Pi$, then $f\in C^3(\Pi)$ and
(\ref{E-18-R}) holds in~$\Pi$.

(ii) If (\ref{E-02-R}) holds for $x_2=0$ and (\ref{E-18-R}) holds in~$\Pi$,
then (\ref{E-02-R}) holds in $\Pi$.
\end{proposition}

\textbf{Proof}.
(i) If (\ref{E-02-R}) is satisfied in $\Pi$, then, since
$H_{ij}$ are of  class $C^1$,
equations (\ref{E-02-R},b) imply that $p,q$ are of class $C^2$.
Hence $f\in C^3(\Pi)$ and, by commutativity of partial derivatives,
(\ref{E-18-R}) holds in $\Pi$.

(ii) Denote for short $x_1=u,\,x_2=v$. By (\ref{E-02-R}) with $v=0$ and (\ref{E-04-R}),
we conclude that $p,q$ and $f$ satisfy in $\Pi$ the integral equations
\begin{subequations}
\begin{eqnarray}
 \label{E-31b}
\begin{array}{c}
  p(u,v) = p(0,0)+\int\nolimits_{0}^v H_{12}(u,\eta,f(u,\eta),p(u,\eta),q(u,\eta))\,d\eta\\
         +\int\nolimits_{0}^u H_{11}(\xi,0,f(u,\eta),p(\xi,0),q(\xi,0))\,d\xi,\\
\end{array}\\
 \label{E-31c}
\begin{array}{c}
  q(u,v) = q(0,0)+\int_{0}^v H_{22}(u,\eta,f(u,\eta),p(u,\eta),q(u,\eta))\,d\eta\\
         +\int\nolimits_{0}^u H_{12}(\xi,0,f(u,\eta),p(\xi,0),q(\xi,0))\,d\xi,\\
\end{array}\\
\label{E-31d}
\begin{array}{c}
  f(u,v) = f(0,0)+\int_{0}^v q(u,\eta)\,d\eta+\int\nolimits_0^u p(\xi,0)\,d\xi
\end{array}\qquad\qquad
\end{eqnarray}
\end{subequations}
(for short we omit the variables $k_1$ and $k_2$ in $H_{ij}$).
By conditions imposed on $p$, $q$, $f$ and $k_i$,
one may differentiate by $u$ the first integrand in (\ref{E-31b}).
Using (\ref{E-04-R}) and
 $\frac{d}{d\,v}H_{11} =H_{11,2} +H_{11,p}\,p_{v}(x,y) +H_{11,q}\,q_{v}(u,v)
 +H_{11,f}\,f_{v}(u,v)$,
we~get
\begin{eqnarray}\label{E-32}
\begin{array}{c}
 \hskip-10mm  p_{u}(u,v) = H_{11}(u,0,f(u,0),p(u,0),q(u,0))\\
 \hskip-10mm +\int\nolimits_{0}^v\big[H_{12,u}(u,\eta,f(u,\eta),p(u,\eta),q(u,\eta))\\
 \hskip-10mm +H_{12,p}(u,\eta,f(u,\eta),p(u,\eta),q(u,\eta))p_{u}(u,\eta)\\
 \hskip-10mm +H_{12,q}(u,\eta,f(u,\eta),p(u,\eta),q(u,\eta))q_{u}(u,\eta)
   +H_{12,f}f_{u}(u,\eta)\big]\/d\eta\\
 \hskip-10mm =H_{11}(u,0,f(u,0),p(u,0),q(u,0))\\
 \hskip-10mm
    +\int\nolimits_{0}^v\frac{d}{d\eta}H_{11}(u,\eta,f(u,\eta),p(u,\eta),q(u,\eta))\,d\eta\\
 \hskip-3mm +\int\nolimits_{0}^v\big[H_{12,p}(p_{u}{-}H_{11})
    {+}H_{12,q}(q_{u}{-}H_{12}){+}H_{12,f}(f_{u}-p){+}H_{12,u}{+}H_{12,p}H_{11}\\
 \hskip-3mm +H_{12,q}H_{12}+H_{12,f}p-H_{11,v}{-}H_{11,p}H_{12} -H_{11,q}H_{22}
    -H_{11,f}q\big]\,d\eta.
\end{array}
 \end{eqnarray}
Define the functions
 $\Theta_1=p_{u}{-}H_{11}(u,v,f,p,q),\
  \Theta_2=q_{u}{-}H_{12}(u,v, f,p,q)$,
 and
 $\Theta_3=f_{u}{-}p.$
By (\ref{E-18-R})$_1$, from (\ref{E-32}) it follows
\begin{subequations}
\begin{equation}\label{E-34}
\begin{array}{ccc}
 \Theta_1(u,v)\eq\int\nolimits_{0}^v [H_{12,p}(u,\eta,f(u,\eta),p(u,\eta),q(u,\eta))
 \,\Theta_1(u,\eta)\\
 \plus H_{12,q}(u,\eta,f(u,\eta),p(u,\eta),q(u,\eta))\,\Theta_2(u,\eta)\\
 \plus H_{12,f}(u,\eta,f(u,\eta),p(u,\eta),q(u,\eta))\,\Theta_3(u,\eta)]\,d\eta.
\end{array}
\end{equation}
Similarly, differentiating (\ref{E-31c}) by $u$, and using (\ref{E-04-R}) and (\ref{E-18-R})$_2$, we obtain
\begin{eqnarray}\label{E-35}
\begin{array}{ccc}
 \Theta_2(u,v)\eq\int\nolimits_{0}^v[H_{22,p}(u,\eta,f(u,\eta),p(u,\eta),q(u,\eta))\,
 \Theta_1(u,\eta)\\
 \plus H_{22,q}(u,\eta,f(u,\eta),p(u,\eta),q(u,\eta))\,\Theta_2(u,\eta)\\
\plus H_{22,f}(u,\eta,f(u,\eta),p(u,\eta),q(u,\eta))\,\Theta_3(u,\eta) ]\,d\eta.
\end{array}
\end{eqnarray}
Differentiating (\ref{E-31d}) by $u$ and using the first equation in (\ref{E-04-R}),
we obtain
\[
\begin{array}{c}
  \hskip-60mm f_{u}(u,v)=p(u,0)+\int\nolimits_0^v p_{v}(u,\eta)\,d\eta\\
  +\int\nolimits_0^v[q_{u}(u,\eta)
 -H_{12}(u,\eta,f(u,\eta),p(u,\eta),q(u,\eta))]\,d\eta.
\end{array}
\]
Hence the following equation is satisfied:
\begin{equation}\label{E-36}
 \Theta_3(u,v)=\int\nolimits_0^v\Theta_2(u,\eta)\,d\eta.
\end{equation}
\end{subequations}
For each $u\in[-a_1,a_1]$ the system of integral equations (\ref{E-34}-c) is equivalent
to Cauchy problem for linear homogeneous ODE's with
initial conditions $\Theta_{i|\,v=0}=0\;(i=1,2,3)$.
Hence $\Theta_1\equiv\Theta_2\equiv\Theta_3\equiv 0$, and (\ref{E-02-R}) are satisfied.
\hfill$\square$

\vskip1.0mm
First, we will recover the graph $M^2$ of $f:\Pi\to[-a_3,a_3]$ infinitesimally along $\gamma$,
i.e., to solve (\ref{E-pqfa1}) for $f_0$ and $p_0,\,q_0$.
 Define the quantity
\begin{equation}\label{dfnbrhoA}
 \bar k_0=\max\{\|\bar k_{1}-\lambda_{|\,\gamma}\|_{\,\gamma},\
 \|\bar k_{2}-\lambda_{|\,\gamma}\|_{\,\gamma}\}.
\end{equation}

\begin{lemma}\label{L-04}
Let $\alpha_0=\alpha_{\,|\,\gamma}\in C^{1}(\gamma)$.
Then for any $r\in(0,a_3]$ there is $\Delta\in(0,r]$ such that for
$\bar k_1,\bar k_2\in C^{0}(\gamma)$ satisfying $\|\bar k_i-\lambda_{|\,\gamma}\|_{\,\gamma}<\Delta\;(i=1,2)$,
the Cauchy problem
\begin{eqnarray}\label{E-pqfa1}
\nonumber
 \ d f_{0}/d x_1\eq p_0,\\
\nonumber
 d p_{0}/d x_1 \eq H_{0,11}^{(1)}\bar k_{1} +H_{0,11}^{(2)}\bar k_{2} -L_0,\\
\nonumber
 d q_{0}/d x_1 \eq H_{0,12}^{(1)}\bar k_{1} +H_{0,12}^{(2)}\bar k_{2} -M_0,\\
 f_0(0)\eq p_0(0) = q_0(0)=0
\end{eqnarray}
has on $\gamma$ a unique $C^{1}$-regular solution
$(f_0,p_0,q_0)$ satisfying $\|(f_0,p_0,q_0)\|_{\,\gamma} < r$.
\end{lemma}

\textbf{Proof}. Denote $\Delta_i=\bar k_i-\lambda_{\,|\,\gamma}$ for $i=1,2$.
Substituting $\bar k_i=\lambda_{\,|\,\gamma}+\Delta_i$ into
(first two equations of) (\ref{E-pqfa1}) and using (\ref{E-tumb4}) gives us
along $\gamma$
\begin{eqnarray}\label{E-tumb3b}
\nonumber
 {d f_{0}}/{d x_1}\eq p_0,\\
\nonumber
 {d p_{0}}/{d x_1}\eq\delta_0\lambda_{\,|\,\gamma} E_0+L_0
 +H_{0,11}^{(1)}\,\Delta_1 +H_{0,11}^{(2)}\,\Delta_2,\\
\nonumber
 {d q_{0}}/{d x_1}\eq\delta_0\lambda_{\,|\,\gamma}\alpha_0F_0+M_0
 +H_{0,12}^{(1)}\,\Delta_1 +H_{0,12}^{(2)}\,\Delta_2,\\
 \label{E-ODEgamma-4}
 f_0(0)\eq p_0(0) = q_0(0)=0
\end{eqnarray}
with
 $E_0=E_{|\,\gamma}, F_0=F_{|\,\gamma}, G_0=G_{|\,\gamma}$,
 and
 $\delta_0=\big(\frac{E_0 G_0-F_0^2}{\det\hat g_{\,|\gamma}}\big)^{1/2}$.
Due to Example~\ref{Ex-Pi},
$\delta_0\lambda_{\,|\,\gamma} E_0+L_0=\delta_0\lambda_{\,|\,\gamma}\alpha_0F_0+M_0=0$
on $\Pi$, see (\ref{E-tumb4}).
Hence $f_0\,{=}\,p_0\,{=}\,q_0\equiv0$ is the solution to (\ref{E-tumb3b}) with $\Delta_i\equiv0$.
Let us take $r\in(0,a_3]$.
We claim that if $\bar k_{i}$ are close enough to $\lambda_{|\,\gamma}$,
then the Cauchy problem (\ref{E-tumb3b}) has on $\gamma$ a unique smooth solution
satisfying $\|(f_0,p_0,q_0)\|_\gamma < r$.
 In aim to apply Proposition~\ref{P-twoODEs}
to (\ref{E-tumb3b}) and to the same system with $\Delta_i\equiv 0$
($f_0\,{=}\,p_0\,{=}\,q_0\equiv 0$ in the last case), denote by
\begin{eqnarray*}
 && P=\left(\begin{array}{c}
           \delta_0\lambda_{\,|\,\gamma} E_0 +L_0 \\
           \delta_0\lambda_{\,|\,\gamma}\alpha_0 F_0+M_0 \\
           p_0\\
          \end{array}\right),\quad
    Q=P+\left(\begin{array}{c}
    H_{0,11}^{(1)}\,\Delta_1+H_{0,11}^{(2)}\,\Delta_2 \\
    H_{0,12}^{(1)}\,\Delta_1+H_{0,12}^{(2)}\,\Delta_2 \\
     0\\
 \end{array}\right).
\end{eqnarray*}
Since the functions $P$ and $Q$ have continuous partial derivatives
w.r. to $f_0,p_0,q_0$, they satisfy in
$\bar\Omega_r=\{|x_1|\le a_1,\ |f_0|\le r,\ |p_0|\le r,\ |q_0|\le r\}$
the Lipschitz condition (for $f_0,p_0,q_0$) with some $\bar L=\bar L(r)>0$.
 By (\ref{E-alpha-3})$_1$,
$\min_{\bar\Omega_r}(\alpha^2 E+2\,\alpha F+G)>0$,
and by (\ref{E-HijkR}), there is $C(r)>0$
such that $|H^{(k)}_{0,ij}|\le\,C(r)$ on $\bar\Omega_r$.
Hence,
\begin{equation}\label{estQmP}
 \| Q - P\|_{\bar\Omega_r}\le 2\,\bar k_0\,C(r),
 \end{equation}
where $\bar k_0$ is defined by (\ref{dfnbrhoA}). Assume that
\begin{equation}\label{E-rhoA-ineq}
 \bar k_0 < \Delta(r):=
 \min\big\{\frac{r}{4\,a_1\,C(r)}\,e^{-\bar L(r)\,a_1},\ r\big\}.
\end{equation}
By the theory of ODE's, there is a maximal interval
$-\eps_1\le t\le\eps_2\;(\eps_i\in(0,a_1])$, in which
(\ref{E-tumb3b}) admits a unique solution with the property
\begin{equation}\label{prop}
 |f_0(x_1)|\le r/2,\quad |p_0(x_1)|\le r/2,\quad |q_0(x_1)|\le r/2.
\end{equation}
On the other hand, in view of (\ref{estQmP}), (\ref{E-rhoA-ineq})
and Proposition~\ref{P-twoODEs}, this solution
satisfies in $[-\eps_1,\eps_2]$ strong inequalities
$|p_0(x_1)|< r/2$, $|q_0(x_1)|<r/2$ and $|f_0(x_1)|< r/2$.
If $\eps_1< a_1$ or $\eps_2< a_1$, due to the theory of ODE's
the solution can be extended on a larger interval
with the property~(\ref{prop}). Hence, $\eps_1=\eps_2=a_1$.
\hfill$\square$

\vskip1.0mm
\textbf{Proof of Theorem~\ref{T-main0}}.
Since $l$ is transversal to $\gamma$, one may assume
$l=\alpha(x_1,x_2)\bar e_1 + \bar e_2$ for some function
$\alpha$ of class $C^{2}$ in a neighborhood of $\gamma$ in $\Pi$.
Since $l$ is not orthogonal to~$\gamma$, we have $\alpha\,\bar g_{11}+\bar g_{12}\ne 0$
along~$\gamma$.
There is $r_1\in(0,a_3]$ such that $\alpha\,\hat g_{11}+\hat g_{12}\ne 0$ over $\gamma$
for $|f|\le r_1$.
Define the functions $E,F,G$ and $L,M,N$ by (\ref{E-EFG},c) in what follows.
In view of $ |(\alpha E+F)-(\alpha\,\hat g_{11}+\hat g_{12})|\le \hat g_{33}|(\alpha p+q)p|$,
see (\ref{E-alpha-3})$_2$, we get $\alpha E+F\ne0$ on the set
$\{(x_1,0,f,p,q): |x_1|\le a_1, \|(f,p,q)\|_\infty\le r\}$ for some $r\in(0,r_1]$.

Restricting (\ref{E-02-R}) on $\gamma$ and denoting $\alpha_0=\alpha_{|\,\gamma}$, $M_0=M_{|\,\gamma},L_0=L_{|\,\gamma}$ and $H_{0,ij}^{(k)}:=H_{ij\,|\,\gamma}^{(k)}$,
yields the system (\ref{E-pqfa1}) for the functions $f_{0}$ and $p_{0}$, $q_{0}$.
 By~Lemma~\ref{L-04}, there exist $\Delta\in(0,r]$ such that if
$\|\bar k_{i}-\lambda_{|\,\gamma}\|_{\,\gamma}<\Delta\;(i=1,2)$,
the Cauchy problem (\ref{E-pqfa1}) has on $[-a_1,a_1]$ a unique solution $(f_0,p_0,q_0)$ of class $C^1$, satisfying $\|(f_0,p_0,q_0)\|_{\,\gamma}<r$.
 By Propositions~\ref{T-compR} and \ref{T-01}, the \textbf{Problem~\ref{Prob-1}} is reduced
 to the Cauchy problem, see (\ref{E-04-R}) and (\ref{E-3.2-R}),
\begin{eqnarray}\label{E-system-5PDEs}
\nonumber
 &&\hskip-8mm f_{x_2}=q(x_1,x_2),\\
\nonumber
 &&\hskip-8mm p_{x_2}= H_{,12}^{(1)}\,k_1 +H_{,12}^{(2)}\,k_2 -M_1,\\
\nonumber
 &&\hskip-8mm q_{x_2}= H_{,22}^{(1)}\,k_1 +H_{,22}^{(2)}\,k_2 -N_1,\\
\nonumber
 &&\hskip-8mm k_{1,x_2}-\frac{\alpha F+G}{\alpha E+F}\,k_{1,x_1}
 =\Psi_1(x_1,x_2, f,p,q, k_1, k_2),\\
 &&\hskip-8mm k_{2,x_2}+\alpha\,k_{2,x_1}= \Psi_2(x_1,x_2, f,p,q,\, k_1,\, k_2),
\end{eqnarray}
with the initial conditions for the functions $f,p,q$ and $k_1,k_2$,
\begin{subequations}
\begin{eqnarray}\label{E-bc0-BB2}
 f(\cdot,0)\eq f_{0},\qquad
 p(\cdot,0)= p_{0},\qquad
 q(\cdot,0)= q_{0},\\
\label{initbrk}
 k_i(\cdot,0)\eq\bar k_i,\quad (i=1,2).
\end{eqnarray}
\end{subequations}
Notice that there exists $\rho\in(0,a_2]$ such that
$\alpha E+F\ne0$ if $|x_1|\le a_1, |x_2|\le\rho, \|(f,p,q)\|_\infty\le r$,
and (\ref{E-system-5PDEs}) is non-singular.
Since $\alpha\in C^{2}$, from the definition of $\Psi_i$ (see the proof of Proposition~\ref{T-compR}) it follows that the functions
$\Psi_i(x_1,x_2, f,p,q, k_1, k_2)\in C^1(\Omega_{\rho,r})$, where
\[
 \Omega_{\rho,r}:=\{(x_1,x_2, f,p,q, k_1,k_2):\ |x_1|\le a_1,\ |x_2|\le\rho,\
 \|(f,p,q)\|_\infty\le r\}.
\]
The normal curvature $\lambda$ of a $C^3$-surface $\Pi$ is $C^1$-regular.
From above it follows that functions in right hand side of (\ref{E-system-5PDEs}) belong to class $C^{1}(\Omega_{\rho,r})$.
The system (\ref{E-system-5PDEs}) for the functions $f,p,q,k_1,k_2$ is hyperbolic,
and it has a diagonal form in its main part
(containing the derivatives of unknown functions).
First three families of characteristics of (\ref{E-system-5PDEs}) are lines $\{x_1=c\}$
and the last two families of characteristics are integral curves of ODE's
$\,\frac{dx_1}{dx_2}=-\frac{\alpha F+G}{\alpha E+F}$ and
$\frac{dx_1}{dx_2}=\alpha$. Denote
$K:=\max\{\|\frac{\alpha F+G}{\alpha E+F}\|_{\tilde\Omega},\,
\|\alpha\|_\Pi\}$, where $\tilde\Omega$ is the projection of $\Omega_{\rho,r}$ onto the space of variables $f,p,q$.
 By Theorem~A (and remark after it, with $c_i>\|\bar k_i\|_{\gamma}$ for $k_i$) there is $\eps\in(0,\rho]$ such that the Cauchy problem (\ref{E-system-5PDEs}), (\ref{E-bc0-BB2},b) admits a unique solution
$\{f,p,q, k_1, k_2\}\in C^{1}(\Pi_{K,\eps})$ with $\|(f,p,q)\|_{\Pi_{K,\eps}}\le r$.
By Proposition~\ref{T-01}, equations (\ref{E-02-R},b)
are valid in $\Pi_{K,\eps}$ and $f\in C^3(\Pi_{K,\eps})$.
Furthermore, by Proposition~\ref{T-compR}, the
functions $k_i\;(i=1,2)$ are the
principal curvatures of the graph $M^2:\,x_3=f(x_1,x_2)$, and, in view of
(\ref{initbrk}), the conditions $k_{i|\,\gamma}=\bar k_i\;(i=1,2)$ are satisfied.
Thus, the surface $M^2$ represents a solution
of \textbf{Problem~\ref{Prob-1}}.

Suppose that a $C^3$-regular surface-graph $M^2:\,x_3=f(x_1,x_2)$
over $\Pi_{K,\eps}$ is a solution of
\textbf{Problem~\ref{Prob-1}} (with principal curvatures $k_i,\;i=1,2$)
satisfying the condition
$\|(f, f_{x_1},f_{x_2})\|_{\Pi_{K,\eps}}\le r$,
where $K$, $\eps$ and $r$ have been chosen above.
By Propositions~\ref{T-compR} and \ref{T-01},
$(f, f_{x_1},f_{x_2},\,k_1,\,k_2)$ is a solution to
(\ref{E-system-5PDEs}) and (\ref{E-bc0-BB2},b),
in which $(f_0,p_0,q_0)$ is a solution to (\ref{E-pqfa1}).
Since solutions of these Cauchy problems are unique,
the solution of \textbf{Problem~\ref{Prob-1}} is
also unique in the class of functions under consideration.
\hfill$\square$

\vskip1.0mm
\textbf{Proof of Corollary~\ref{T-main3}}.
Let $M^2$ in $\R^3$ (with cartesian coordinates) be a surface-graph of $f\in C^{2}(\Pi)$ defined on $\Pi=\{|x|\le a_1,\,|y|\le a_2,\,z=0\}$.
By Lemma~\ref{L-EFG-1} or directly, we find that
 $\,n= \frac{1}{\sqrt{1+p^2+q^2}}[-p,\,-q,\,1]$
is the unit normal to $M^2$, and the 1-st and the 2-nd fundamental forms of $M^2$ are
 $E= 1+p^2,\,F= p\,q,\,G= 1+ q^2$, and
 $L= \frac{f_{x x}}{\sqrt{1{+}p^2{+}q^2}},\
  M= \frac{f_{x y}}{\sqrt{1{+}p^2{+}q^2}},\
  N= \frac{f_{y y}}{\sqrt{1{+}p^2{+}q^2}}$.
From Proposition~\ref{T-compR} it follows that
$k_1,k_2\in C^{0}(\Pi)$ are the principal curvatures and the vector $l=(\alpha(x,y),1)$ (where $\alpha\in C^{1}(\Pi)$) is the projection onto $\Pi$ of $k_1$-principal direction of $M^2$ if and only~if
\begin{subequations}
\begin{eqnarray}\label{E-02}
 p_{x} \eq H_{11}(x,y,p,q,k_1,k_2),\quad
 q_{x} = H_{12}(x,y,p,q,k_1,k_2),\\
\label{E-04}
 p_{y} \eq H_{21}(x,y,p,q,k_1,k_2),\quad
 q_{y} =H_{22}(x,y,p,q,k_1,k_2).
\end{eqnarray}
\end{subequations}
For $\R^3(x,y,z)$, the system (\ref{E-04-R})\,--\,(\ref{E-3.2-R}) with $\lambda\,{=}\,0$
has a form (see also (\ref{E-HijkR-R3}))
\begin{eqnarray}\label{E-system-4PDEs}
\nonumber
 && p_{y} = H_{,12}^{(1)}\,k_1 +H_{,12}^{(2)}\,k_2,\\
\nonumber
 && q_{y}= H_{,22}^{(1)}\,k_1 +H_{,22}^{(2)}\,k_2,\\
\nonumber
 && k_{1,y}-\frac{\alpha\,p\,q+q^2+1}{\alpha(p^{2}{+}1)+p\,q}\,k_{1,x}
 +(c_{1}+c_3 k_1)(k_2-k_1)= 0,\\
 && k_{2,y}+\alpha\,k_{2,x}+c_{2}(k_{2}-k_{1})= 0,
\end{eqnarray}
where
 $c_{1}=\frac{(p^{2}{+}q^{2}{+}1)(\alpha\,{\alpha}_{,x}+{\alpha}_{,y})}
 {\delta_1(\alpha(p^2+1)+p\,q)}$,
 $c_{2}=\frac{(\alpha({p}^2{+}1){+}p\,q){\alpha}_{,y}
 {-}(\alpha p\,q{+}{q}^{2}{+}1){{\alpha}_{,x}}}{\delta_1}$,
 $c_3={\frac{\sqrt{p^{2}+q^{2}+1}(\alpha\,q-p)}{\alpha(p^{2}+1)+p\,q}}$.
If $f\in C^3(\Pi)$ and $k_1,k_2\in C^1(\Pi)$, then (\ref{E-system-4PDEs})$_{3,4}$
are compatibility conditions for (\ref{E-02},b), i.e.,
$(p_{x})_{y}=(p_{y})_{x}$, $(q_{x})_{y}=(q_{y})_{x}$.

Similarly to the proof of Theorem~\ref{T-main0}, one may show that
there are $r,K>0$ and $\Delta\in(0,r)$ such that
if $\|\bar k_{i}\|_\gamma<\Delta$, then for some $\eps\in(0,a_2)$
the Cauchy problem (\ref{E-system-4PDEs}) with
 $p(\cdot,0)= p_{0},\ q(\cdot,0)= q_{0},\ k_i(\cdot,0)=\bar k_i,\ (i=1,2)$
is non-singular and admits a unique solution $p,q,k_1,k_2$ of class $C^{1}(\Pi_{K,\eps})$
with $\|(p,q,k_1,k_2)\|_{\Pi_{K,\eps}}\le r$.
Moreover, there is a unique function $f$ of class $C^{2}(\Pi_{K,\eps})$
such that $f_{x}= p,\ f_{y}= q$ and $f(0,0)=df(0,0)=0$.
As in the proof of Theorem~\ref{T-main0}, one may show
that $f\in C^3(\Pi_{K,\eps})$ and the surface-graph $M^2:\,z=f(x,y)$ represents
a unique solution of \textbf{Problem~\ref{Prob-1}} in the class of functions
from the formulation of the corollary.
\hfill$\square$

\subsection{Proof of Theorem~\ref{T-RPC2} and Proposition~\ref{T-RPC1}}
\label{sec:Th2-3}

\vskip1.0mm
A surface $M^{2}$ in $\R^3(k)$ is called \textit{parallel curved} (PC) if
there is a totally umbilical (or totally geodesic) surface $\beta\subset\R^3(k)$ such that at each point $x\in M^2$ there is a principal direction tangent to $\beta_d$ (a parallel surface to $\beta$ on the distance $d$).
Surfaces of revolution and cylinders in Euclidean space $\R^3$ provide examples of PC surfaces (see \cite{Ando2004} and~\cite{Ando2007}).
 We study PC surfaces in space forms $\R^3(k)$ (see also Section~\ref{sec:warped2}) in relation to \textbf{Problems~\ref{Prob-1}}
and \textbf{\ref{Prob-2}}.

A PC surface $M^2$ in $\R^3(k)$ (a Riemannian 3-space of constant curvature $k$)
can be recovered locally by exact procedure.
 We will prove the claim for PC surfaces in~$\R^3$.
 In Proposition~\ref{T-RPC1} we solve \textbf{Problem~\ref{Prob-2}}
using twice the reconstruction a plane curve by its curvature.

\vskip1.0mm
\textbf{Proof of Proposition~\ref{T-RPC1}} will be divided into three steps.
\newline
1. We apply geometrical construction of Proposition~\ref{L-RPC-1a} when $\beta\subset\R^3$
is a plane with the normal $\alpha\,e_2-e_1$.
In this case, $l=\alpha\,e_1+e_2$ is a constant vector field
(its trajectories are parallel lines). The planes parallel to $\beta$
intersect $M^2$ transversally by curvature lines, Fig.~\ref{F-00y}(a).

Let $\pi_1:M^2\to\beta_0$ be the orthogonal projection onto the plane $\beta_0=\{x-\alpha\,y=0\}$.
 Let $\gamma_0(t)\subset\beta_0$ be a
$k_1$-curvature line on $M^2$ through the origin $O(0,0,0)$.
The normals to $\gamma_0$ (lines in $\beta_0$)
and parallel curves (of constant distance) to $\gamma_0$ form a {\em semi-geodesic net} on $\beta_0$ on a neighborhood of $\gamma_0$.
Let us parameterize~$\gamma_0$
\begin{equation}\label{E-v0-gamma0}
 v=v_0(u),\qquad v_0(0)=v'_0(0)=0\quad(|u|\le \eps\ \mbox{ for some } \eps>0)
\end{equation}
in $(u,v)$-coordinates of $\beta_0$.
 Notice that $\tilde e_1,\tilde e_2,\tilde e_3$ ($\tilde e_3=(\alpha\,e_1-e_2)/\sqrt{\alpha^2+1}$ is a unit normal to $\beta_0$) is the orthonormal frame of $\R^3$. The curvature of $\gamma_0$ is
\begin{equation}\label{E-ku}
 \tilde k(u)={v_0''(u)}/{(1+{v_0'}^2(u))^{3/2}}.
\end{equation}
From (\ref{E-v0-gamma0}) and (\ref{E-ku}) it follows the inequality
$|v_0'(u)|\le\bar k\int_{0}^u{[1+{v_0'}^2(s)]^{3/2}}\,ds $.
By Lemma~\ref{appl}, if $\bar k\le\frac{\sqrt{\alpha^2+1}}{\alpha\,a_1}$,
then $|v_0'(u)|\le \bar k\,u/{(1-(\bar k\,u)^2)^{1/2}}$
(a'priori estimate) for $|u|\le \frac{a_1\,\alpha}{\sqrt{\alpha^2+1}}$.
If we take $\bar k\le\frac{\sqrt{\alpha^2+1}}{\sqrt{2}\,\alpha\,a_1}$,
then a unique solution $v_0(u)$ to (\ref{E-v0-gamma0}), (\ref{E-ku}),
defined for $|u|\le \frac{a_1\alpha}{\sqrt{\alpha^2+1}}$, satisfies the inequalities
\begin{equation}\label{estv0}
 |v_0|\le \sqrt{2}\,\frac{a_1^2\,\bar k\,\alpha^2}{\alpha^2+1},\quad
 |v_0'|\le \sqrt{2}\,\frac{a_1\bar k\,\alpha}{\sqrt{\alpha^2+1}},\quad
 |v_0''|\le\sqrt8\,\bar k.
\end{equation}
Since $\tilde k(u)$ is $C^1$-regular, from (\ref{E-ku}) it follows that
$v_0(u)$ is $C^3$-regular. Hence, the curve
 $\gamma_\rho(u)=[u-\frac{v_0'(u)\,\rho}{[1+{v_0'}^{2}(u)]^{1/2}},\
 v_0(u)+\frac{\rho}{[1+{v_0'}^{2}(u)]^{1/2}}]$,
 (on the distance $\rho$ to $\gamma_0$ in the plane $\beta_0$)
is $C^{\,2}$-regular for small enough $\rho$, its curvature is
$k_\rho(u)=\frac{-\tilde k(u)}{1-\,\tilde k(u)\rho}$, see~\cite{Top}.

2. Assume now that $\rho=\rho(h)$ (the function of $h$) and translate $\gamma_{\rho(h)}$
 on the height $h$ in the normal direction to $\beta_0$.
 Then, the obtain $(u,h)$-parametrizati\-on of $M^2$ near $\gamma_0$ with $\beta_0$-level curves~$\gamma_\rho$,
 \begin{equation}\label{E-mpar}
 U=u-\frac{v_0'(u)\,\rho(h)}{(1{+}{v_0'}^{2})^{1/2}},\quad
 V=v_0(u)+\frac{\rho(h)}{(1{+}{v_0'}^{2})^{1/2}},\quad
 W=h.
\end{equation}
Using the equation for principal curvatures,
$(EG-F^2)k^2-(EN+GL-2FM)k+(LN-M^2)=0$, since $F=M=0$, we find
 \begin{equation}\label{E-k12}
 k_1(u,h) = \frac{\tilde k(u)}{1-\tilde k(u)\rho(h)}\cdot
 \frac{1}{({1+{\rho'}^2(h)})^{1/2}}, \ \
 k_2(u,h) = \frac{\rho''(h)}{({1+{\rho'}^2(h)})^{3/2}}.
\end{equation}
Notice that $k_1(u,h)$ is the curvature of $\gamma_\rho$ multiplied by $\cos\varphi$, see (\ref{E-k12})$_1$, where $\varphi$ is the angle between $\beta_0$ and the tangent plane to $M^2$ through the intersection point.
The 2-nd principal curvature, $k_2(u,h)$ (does not depend on $u$) is simply the curvature of the curve ${\bold r}(h)=[h,\rho(h)]$ in the vertical plane $\beta^\perp=\{\alpha x+y=0\}$ with the $\rho$-axis $Oz$.
 One may recover a $C^2$-regular function $\rho(h)$ for small enough $h$ from (\ref{E-k12})$_2$, solving the BVP
\begin{equation}\label{E-rhok2}
  {\rho''(h)}{({1+{\rho'}^2(h)})^{3/2}}=\bar k_2(h\sqrt{\alpha^2+1}),\qquad
  \rho(0)=\rho'(0)=0.
\end{equation}
If $\bar k\le\frac{\sqrt{\alpha^2+1}}{\sqrt{2}a_1}$, then by Lemma~\ref{appl},
a unique solution $\rho(h)$ to Cauchy problem (\ref{E-rhok2}) exists for
$|h|\le \frac{a_1}{\sqrt{\alpha^2+1}}$, and satisfies the inequalities
 \begin{equation}
\label{estrh1}
 |\rho|\le \sqrt{2}\frac{a_1^2\,\bar k}{\alpha^2+1},\quad
 |\rho'|\le \sqrt{2}\frac{a_1\bar k}{\sqrt{\alpha^2+1}},\quad
 |\rho''|\le\sqrt8\,\bar k.
\end{equation}
One may recover $\tilde k(u)$ from (\ref{E-k12})$_1$ for $\rho=0$
with known $k_1(u,0)$.
Finally, $v_0(u)$ is a unique $C^2$-regular solution to (\ref{E-ku})
with initial values~(\ref{E-v0-gamma0}).
It was shown that \textit{if $\bar k\le \frac{\sqrt{\alpha^2+1}}{\sqrt{2}\,a_1}\min\{1, 1/\alpha\}$, then the parametrization (\ref{E-mpar}) defines a regular surface over the rectangle of parameters $(u,h)$}
$$
 \tilde\Pi(a_1)=\big\{(u,h)\in\R^2:\
 |u|\le{a_1\alpha}/{\sqrt{\alpha^2+1}},\ |h|\le{a_1}/{\sqrt{\alpha^2+1}}\big\}.
$$
Equations (\ref{E-mpar}) of $M^2$ in cartesian coordinates $(x,y,z)$ take a form
\begin{equation}\label{E-mpar2}
 X=\frac{h+\alpha\,U(u,h)}{\sqrt{\alpha^2+1}},\quad
 Y=\frac{U(u,h)-\alpha\,h}{\sqrt{\alpha^2+1}},\quad
 Z=V(u,h).
\end{equation}
 In view of (\ref{E-mpar}), (\ref{E-mpar2}),
 $C^3$-regularity of $v_0(u)$ and $C^2$-regularity of $\rho(h)$,
 the surface $M^2$ is $C^2$-regular.

3. Along $\gamma_1$ (i.e., $\rho=0$) we get
$[X_{,h}, Y_{,h}, Z_{,h}]_{|(u,0)}=
 [\frac{1}{\sqrt{\alpha^2+1}},\ \frac{-\alpha}{\sqrt{\alpha^2+1}},\ 0]$,
 and
 ${[X_{,u}, Y_{,u}, Z_{,u}]}_{|(u,0)}=
 [\frac{\alpha}{\sqrt{\alpha^2+1}},\ \frac{1}
 {\sqrt{\alpha^2+1}},\ v_0'(u)].$
Hence
  $\det\Big(\begin{array}{cc}
        X_{,h}(u,0) & Y_{,h}(u,0) \\
        X_{,u}(u,0) & Y_{,u}(u,0) \\
      \end{array}\Big)$ $=1>0$,
and $M^2$ regularly projects onto a neighborhood of $\gamma_1$ in $xy$-plane.
It is not difficult to see that $M^2$ is $C^2$-regular.
Let us show that for any $a\in(0,a_1)$ there is $\delta>0$
such that if $\bar k<\delta$, the surface $M^2$ regularly projects onto $\Pi(a)$.
Based on (\ref{E-mpar})$_{1}$, we write (\ref{E-mpar2})$_{1,2}$ in the equivalent form
 $u=\frac{\alpha\,X+Y}{\sqrt{\alpha^2+1}}+\frac{v_0'(u)\,\rho(h)}
 {(1+{v'_0}^2)^{1/2}},\ h=\frac{X-\alpha\,Y}{\sqrt{\alpha^2+1}}$.

Consider the mapping
 $T:(u,h)\to\big[\frac{\alpha\,x+y}{\sqrt{\alpha^2+1}}+\frac{v_0'(u)\,\rho(h)}
 {(1+{v'_0}^2)^{1/2}},\; \frac{x-\alpha\,y}{\sqrt{\alpha^2+1}}\big]$
associated with above system. We~will show that for any $(x,y)\in\Pi(a)$ the system
\[
 u=\frac{\alpha\,x+y}{\sqrt{\alpha^2+1}}+\frac{v_0'(u)\,\rho(h)}
 {(1+{v'_0}^2)^{1/2}},\quad h=\frac{x-\alpha\,y}{\sqrt{\alpha^2+1}}.
\]
admits a unique solution $(u,h)$ in $\tilde\Pi(a_1)$.
Notice $\tilde\Pi(a)$ is the image of $\Pi(a)$ under the linear mapping
$u=\frac{\alpha\,x+y}{\sqrt{\alpha^2+1}}$, $h=\frac{x-\alpha\,y}{\sqrt{\alpha^2+1}}$.

Assume that the metric in $\tilde\Pi(a_1)$ is induced by the norm
$\|(u,h)\|_\infty=\max\{|u|,|h|\}$. In order to apply the Banach
fixed point theorem, we will show that for a small enough $\bar k$
the mapping $T$ maps $\tilde\Pi(a_1)$ into itself, and that $T$ is a contraction.
If $\bar k\le\frac{(a_1-a)\,\alpha\,\sqrt{\alpha^2+1}}{\sqrt{2}\,a_1^2}$,
then using (\ref{estrh1})$_1$, we obtain for $(x,y)\in\Pi(a)$
and $(u,h)\in\tilde\Pi(a)$ that ${|x-\alpha\,y|}\le {a_1}$ and
\begin{eqnarray*}
 &&\Big|\frac{\alpha\,x+y}{\sqrt{\alpha^2+1}}+\frac{v_0'(u)\,\rho(h)}
 {(1+{v'_0}^2)^{1/2}}\Big|\le\frac{\alpha\,a}{\sqrt{\alpha^2+1}}+|\rho(h)|
 \le\frac{\alpha\,a_1}{\sqrt{\alpha^2+1}}.
\end{eqnarray*}
Under above condition for $\bar k$, $T$ maps $\tilde\Pi(a_1)$ into itself.
To show that $T:\,\tilde\Pi(a_1)\to\tilde\Pi(a_1)$ is a contraction, we find
the differential of $T$ at $(u,h)\in\tilde\Pi(a_1)$:
\[
 d\,T(u,h)\Big(\begin{array}{c}
 \Delta u \\ \Delta h \end{array}\Big)
 =\Big(\begin{array}{c}
 -\frac{\rho(h)\,v_0''(u)}{(1+{v'_0}^2(u))^{3/2}}\,\Delta u
 +\frac{\rho'(h)\,v_0'(u)}{(1+{v'_0}^2)^{1/2}}\,\Delta h \\
 0
 \end{array}\Big).
\]
Using (\ref{estv0}), (\ref{estrh1}), we obtain the following estimates:
\begin{eqnarray*}
 \Big|\frac{\rho(h)\,v_0''(u)}{(1+{v'_0}^2)^{3/2}}\Big|\le
 \frac{4\sqrt{2}\,a_1^2}{\alpha^2+1}\,\bar k^2,\quad
 \Big|\frac{\rho'(h)\,v_0'(u)}{(1+{v'_0}^2)^{1/2}}\Big|\le
 \frac{\sqrt{2}\,a_1}{\sqrt{\alpha^2+1}}\,\bar k.
\end{eqnarray*}
Hence, for a small enough $\bar k$ the norm $\|dT(u,h)\|_\infty$
is less than $1$ for any $(u,h)\in\tilde\Pi(a_1)$,
that is the mapping $T$ is a contraction in $\tilde\Pi(a_1)$.
Thus, for small enough $\bar k$, the projection ${\rm proj}_{xy}$
(the orthogonal projection onto the $xy$-plane) realizes a bijection between the surface $M^2\cap\mathrm{proj}_{xy}^{-1}(\Pi(a))$ and $\Pi(a)$.
Since the linear operator ${\rm id}-d\,T$ is invertible for any $(u,h)\in\tilde\Pi(a_1)$,
 by the Implicit Function Theorem this projection is regular.\hfill$\square$


\vskip1.5mm
\textbf{Proof of Theorem~\ref{T-RPC2}}.
The condition $Y(u,h)\equiv 0$ determines the intersection curve of $M^2$
(see (\ref{E-mpar2}))
with $xz$-plane of the form $(X(u),0,Z(u))$, where $Z(u):=V(u,\frac{X(u)}{\sqrt{\alpha^2+1}})$.
Now, $X(u)$ and $v_0(u)$ are functions to be found using the boundary values
of the principal curvatures. Denote $w=v_0'(u)$. Consider the system
\begin{eqnarray}\label{eqhu1}
 \hskip-6mm X=\frac{\sqrt{\alpha^2+1}}{\alpha}\,
 \Big(u-\frac{\tilde\rho(X)\,w}{(1+w^2)^{1/2}}\Big),\quad
 w=\int\nolimits_0^u\frac{\bar k_1(X)(1+w^2)^{3/2}\,d\eta}
 {\tilde\rho(X)\bar k_1(X)+\Phi(X)}\,,
\end{eqnarray}
where $\tilde\rho(x):=\rho(\frac{x}{\sqrt{\alpha^2+1}})$,
$\Phi(x):=\tilde\Phi(\frac{x}{\sqrt{\alpha^2+1}})$,
$\tilde\Phi(h)=\frac{1}{(1+{\rho'(h)}^2)^{1/2}}$,
$\bar k_1(x)$ is a continuous in $[-a,a]$ function,
and $X=X(\eta)$, $w=w(\eta)$ in the integrand.

Due to the proof of Proposition~\ref{T-RPC1}, we are looking for the parametric form
(\ref{E-mpar2}) of $M^2$, where $U,V,W$ (the functions of $u,h$)
are given in (\ref{E-mpar}), and $v_0(u)$ satisfies (\ref{E-ku})
with unknown function $\tilde k(u)$.
By Proposition~\ref{T-RPC1}, the principal curvature $k_2$ of
$M^2$ satisfies $k_2(x,0)=\bar k_2(x)$.
It is sufficient to show that for small enough $\bar k_i$ there exists a unique pair
$X(u),\,v_0(u)$ (and hence $M^2$) such that  $k_1(x,0)=\bar k_1(x)$ for the principal curvature $k_1$ of $M^2$.
 From (\ref{E-mpar2}) with $Y=0$ we conclude that the pair $X=(X(u),v_0(u))$
satisfies the equation $\frac{\alpha}{\sqrt{\alpha^2+1}}X=u-\frac{\tilde\rho(X)\,v_0'(u)}{(1+{v'_0}^2)^{1/2}}$
for $u\in I=[-\frac{a\,\alpha}{\sqrt{\alpha^2+1}},\frac{a\,\alpha}{\sqrt{\alpha^2+1}}]$ $(a\in(0,a_1))$.
In order to satisfy $k_1(x,0)=\bar k_1(x)$, we rewrite (\ref{E-k12})$_1$ as
\begin{equation}\label{E-d-u}
 \frac{\tilde k(u)}{1-\tilde k(u)\tilde \rho(X)}\,\Phi(X)=
 \bar k_1(X) \ \Leftrightarrow \
 \tilde k(u) =\frac{\bar k_1(X)}{\tilde\rho(X)\bar k_1(X)+\Phi(X)}.
\end{equation}
Substituting $\tilde k(u)$ of (\ref{E-ku}) into (\ref{E-d-u}),
we get the integral equation
$$
 v_0'(u)=\int\nolimits_0^u\frac{\bar k_1(X(\eta))(1+{v'_0}^2(\eta))^{3/2}} {\tilde\rho(X(\eta))\bar k_1(X(\eta)) +\Phi(X(\eta))}\,d\eta\quad(u\in I).
$$
This leads to the system (\ref{eqhu1}) for the functions $X(u)$ and $w(u):=v_0'$.

We claim that for small enough $\bar k_1$ and $\bar k_2$,
(\ref{eqhu1}) admits a unique solution in the product space $\mathcal{B}:=B(0,a_1)\times C^0(I)$,
where $B(\phi,r)$ is the closed ball of radius $r>0$ in $C^0(I)$
centered at the function $\phi$.
First we will investigate (\ref{eqhu1})$_1$.
 Set $\bar k:=\max\limits_{|x|\le a_1}\{|\bar k_1(x)|,|\bar k_2(x)|\}$.
 By the proof of Proposition~\ref{T-RPC1}, if $\bar k\le\frac{\sqrt{\alpha^2+1}}{\sqrt{2}\,a_1}$,
 the Cauchy problem (\ref{E-rhok2}) admits a unique solution $\rho(h)$ defined
 for $|h|\le a_1/\sqrt{\alpha^2+1}$ and the estimates (\ref{estrh1}) are valid.
 Denote
$$
 F(X,u,w):=\frac{\sqrt{\alpha^2+1}}{\alpha}\,\Big(u-\frac{\tilde\rho(X)\,w}{(1+w^2)^{1/2}}\Big).
$$
 From (\ref{estrh1})$_{1,2}$ we get that if
 $\bar k\le\min\{\frac{(a_1-a)\,\alpha\,\sqrt{\alpha^2+1}}{\sqrt{2}\,a_1^2},
 \frac{\alpha\sqrt{\alpha^2+1}}{2\,a_1}\}$. Hence
 $|F|\le a+\frac{\sqrt{\alpha^2+1}}{\alpha}|\tilde\rho(X)|\le a+\frac{\sqrt{2}\,a_1^2}{\alpha\,\sqrt{\alpha^2+1}}\,\bar k\le a_1$
 for any $(X,u,w)\in[-a_1,a_1]\times I\times\R$, that~is
\begin{eqnarray}\label{propF}
 \forall\;(u,w)\in I\times\R:\quad F(\cdot,u,w):\;[-a_1,a_1]\to[-a_1,a_1],\\
\label{dFdX}
 |\partial_X F|\le\frac{\sqrt{\alpha^2+1}}{\alpha}|\tilde\rho'(X)|\le
 \frac{\sqrt{2}\,a_1}{\alpha\,\sqrt{\alpha^2+1}}\,\bar k\le\frac{1}{\sqrt{2}}.
\end{eqnarray}
Hence $\partial_X(X-F(X,u,w))>0$ for any $(X,u,w)\in [-a_1,a_1]\times I\times\R$.
We conclude that (\ref{eqhu1})$_1$ admits a unique solution $X=\tilde X(u,w)\in[-a_1,a_1]$
for any $(u,w)\in I\times\R$. Since the function $F$ is $C^1$-regular in
$[-a_1,a_1]\times I\times\R$, by the Implicit Function Theorem,
the function $\tilde X$ is $C^1$-regular in $I\times\R$.
Let us substitute $X=\tilde X(u,w)$ into
(\ref{eqhu1})$_2$, which is equivalent to Cauchy problem
\begin{equation}\label{eq-dwdu}
 \frac{d\,w}{d\,u}=\frac{\bar k_1(\tilde X)
 (1+w^{2})^{3/2}}{\tilde\rho(\tilde X)\bar k_1(\tilde X)+\Phi(\tilde X)}\quad
 w(0)=0.
\end{equation}
We will show that (\ref{eq-dwdu}) admits a unique solution $w(u)$ in $C^0(I)$.
By estimates (\ref{estrh1})$_{1,2}$ and definition of function $\Phi$, we have
\begin{equation}\label{estblw}
 |\tilde\rho(\tilde X)\,\bar k_1(\tilde X)+\Phi(\tilde X)|\ge\sqrt{2}
 -\sqrt{2}\frac{a_1^2\,\bar k}{\alpha^2+1}\ge\frac{1}{\sqrt{2}}
\end{equation}
for $\bar k\le\frac{\alpha^2+1}{2\,a_1^2}$. Since $\bar k_1(X)$, $\tilde\rho(X)$
 and $\Phi(X)$ are $C^1$-regular in $[-a_1,a_1]$, the ODE (\ref{eq-dwdu})
 satisfies the conditions required for local existence and
 uniqueness of a solution to Cauchy problem.
 In order to show that the solution to (\ref{eq-dwdu}) does not blow up in
 $I$, we need an a'priori estimate of a solution to (\ref{eqhu1})$_2$ with $X=\tilde X(u,w)$.
 Let $w=w(u)$ be a continuous solution to this equation in $[-c,c]\subseteq I$.
 From (\ref{eqhu1})$_2$ and (\ref{estblw}) it follows
$$
 |w(u)|\le\sqrt{2}\,\bar k\,{\rm sign}(u)\int\nolimits_0^u(1+w^2(\eta))^{3/2}\;d\eta\quad
 (u\in [-c,c]).
$$
By Lemma~\ref{appl}, the estimate
 $|w(u)|\le \frac{\sqrt{2}\,\bar k|u|}{\sqrt{1-2\,\bar k^2 u^2}}\le 1$
is valid in $[-c,c]$, if $\bar k\le\frac{\sqrt{\alpha^2+1}}{2\,\alpha\,a_1}$.
Since the bound for $|w(u)|$ in $[-c,c]$ does not depend on
$c$, the solution $w(u)$ exists and is continuous on $I$.
So, we have proved the claim
($X(u)=\tilde X(u,w(u))$) in the product space $\mathcal{B}$,
moreover, by above a'priori  estimate, the solution belongs
to the smaller space $\tilde\mathcal{B}:=B(0,a_1)\times B(0,1)$.

The desired $C^2$-regular surface $M^2$ is given by (\ref{E-mpar2}), (\ref{E-mpar}),
where $v_0(u)=\int_0^u w(\eta)\,d\eta$.
Using (\ref{E-d-u})$_2$ and Proposition~\ref{T-RPC1}, we obtain that
for any $a\in (0,a_1)$ there exists $\delta>0$ such that if $\bar k<\delta$,
then $M^2$ projects regularly onto $\Pi(a)$ and $M^2$ is a unique surface with
the properties indicated in the theorem.\hfill$\square$

\begin{proposition}\label{T-RPC2-B}
The solution in Theorem~\ref{T-RPC2} can be represented in the form (\ref{E-mpar}), (\ref{E-mpar2}),
where $v_0(u)=\int_0^u w(\eta)\,d\eta$ and w(u)
is the second component of the solution $(X,w)$
to~(\ref{eqhu1}).
Moreover, $(X,w)$ is the limit of the iterated function sequence for the operator
$S$ in $C^0(\R)\times C^0(\R)$
\[
 S: (X,w)\to\Big[\frac{\sqrt{\alpha^2+1}}{\alpha}\,
 \Big(u-\frac{\tilde\rho(X)\,w}{(1+w^2)^{1/2}}\Big),\quad
 \int\nolimits_0^u\frac{\bar k_1(X)(1+w^2)^{3/2}\,d\eta}
               {\tilde\rho(X)\bar k_1(X)+\Phi(X)}\Big]
\]
with the starting point $(\frac{\sqrt{\alpha^2+1}}{\alpha}\,u,\ 0)$.
\end{proposition}

\textbf{Proof}.
Indeed, a fixed point of $S$ is a solution to (\ref{eqhu1}).
Let us prove that for a small enough $\bar k$ the solution $(X(u), w(u))$ to (\ref{eqhu1})
(that determines the surface $M^2$) can be found by an iterative process.
First we will show that for a small enough $\bar k$ the operator $S$
maps $\tilde\mathcal{B}=B(0,a_1)\times B(0,1)$ into itself.
Denote $G(X,w):=\frac{\bar k_1(X)(1+w^2)^{3/2}}{\tilde\rho(X)\bar k_1(X)+\Phi(X)}$.
By (\ref{estrh1})$_{1}$ and (\ref{estblw}), for
$\bar k\le \frac{\sqrt{\alpha^2+1}}{4\,\alpha\,a_1}$ we have
\[
 \big\|\int\nolimits_0^u\,G(X(\eta),w(\eta))\,d\eta\big\|_I
 \le\frac{4\,\alpha\,a_1}{\sqrt{\alpha^2+1}}\;\bar k\le 1,\quad
 {\rm where} \quad (X,w)\in\tilde B.
\]
By this and (\ref{propF}), for a small enough $\bar k$ the operator $S$ maps $\tilde\mathcal{B}$ into itself.

Let us show that $S:\tilde\mathcal{B}\to\tilde\mathcal{B}$ is a contraction
w. r. to some metric on $\tilde\mathcal{B}$ for a small enough $\bar k$.
From (\ref{estrh1})$_{1}$ we obtain
$|\partial_w F|\le\frac{\sqrt{\alpha^2+1}}{\alpha}\,|\tilde\rho(X)|\le
\frac{\sqrt{2}\,a_1^2}{\alpha\,\sqrt{\alpha^2+1}}\,\bar k$,
if $(X,u,w)\in[-a_1,a_1]\times I\times\R$.
This estimate and (\ref{dFdX}) mean that $\partial_XF$ and $\partial_wF$
becomes arbitrary small for a small enough $\bar k$, hence the first component of $S=(S_1,S_2)$
satisfies w.r. to $X,w$ the Lipschitz condition with the Lipschitz constant $L_1\in(0,1)$ for a small enough $\bar k$. Let us compute the differential of the second component $S_2$ (of $S$) at a point $(X(u),w(u))\in\tilde\mathcal{B}$:
\[
 dS_2(X,w)\hskip-1pt\Big(\hskip-3pt\begin{array}{c} \Delta X \\
 \Delta w \end{array}\hskip-4pt\Big)
 \hskip-2pt=\hskip-4pt\int\nolimits_0^u\hskip-3pt
 (\partial_XG(X(\eta),w(\eta))\Delta X(\eta)+\partial_wG(X(\eta),w(\eta))\Delta w(\eta))d\eta
\]
where
$\partial_XG=-\frac{(1+w^2)^{3/2}\bar k_1(X)
 (\tilde\rho'(X)\bar k_1(X)+\bar k_1'(X)\tilde\rho(X)+\Phi'(X))}
 {(\tilde\rho(X)\bar k_1(X)+\Phi(X))^2}+\frac{(1+w^2)^{3/2}\bar k_1'(X)}
 {\tilde\rho(X)\bar k_1(X)+\Phi(X)}$
and
$\;\partial_wG=\frac{3\,\bar k_1(X)w(1+w^2)^{1/2}}
{\tilde\rho(X)\bar k_1(X)+\Phi(X)}$.
Using $\tilde\Phi'=-\frac{\rho''\,\rho'}{(1+{\rho'}^2)^{3/2}}$ and (\ref{estrh1}), (\ref{estblw}),
yields that the function $G$ satisfies the Lipschitz condition w.r. to $X,w$ in $[-a_1,a_1]\times[-1,1]$
with a Lipschitz constant $L_2>0$ for a small enough $\bar k$.
 Notice that $L_2$ is not arbitrary small for a small enough $\bar k$,
 because the expression for $dS_2$ contains the derivative $\bar k_1'$ that is not
 assumed to be small.

 In aim to show that $S$ is a  contraction,
 let us define the following metric in the second component $B(0,1)$
 of the product space $\tilde\mathcal{B}$:
 $d^{\,T}(w_1,w_2)=\max_{u\in I}\,e^{-T|u|}|w_1(u)-w_2(u)|$
 (see \cite{Kr}), where $T>0$ will be  chosen in the sequel.
 Clearly, this metric is equivalent to the original $C^0$-metric
 $d_\infty(w_1,w_2)=\max_{u\in I}\,|w_1(u)-w_2(u)|$ in $B(0,1)$.
 The metric of $\tilde\mathcal{B}$~is
 $$\tilde d^{\,T}((X_1,w_1),(X_2,w_2))=\max\{d_\infty(X_1,X_2),\ d^{\,T}(w_1,w_2)\}.$$
 Let us estimate for $(X_i,w_i)\in\tilde\mathcal{B}\ (i=1,2)$, $u\in I$, $u\ge 0$:
\[
\begin{array}{ccc}
 &&\hskip-4mm
 \big|e^{-Tu}(S_2(X_1,w_1)(u)-S_2(X_2,w_2)(u))\big|
 =\big|\int\nolimits_0^u\,e^{-Tu}(G(X_1(\eta),w_1(\eta))\\
 && -G(X_2(\eta),w_2(\eta)))\,d\eta\big|\le
 L_2\int\nolimits_0^u\,e^{-T(u-\eta)}e^{-T\eta}
 \max\{|X_1(\eta)-X_2(\eta)|,\\
 &&\hskip-4mm|w_1(\eta)-w_2(\eta)|\}\,d\eta\le L_2
 \int\nolimits_0^ue^{-T(u-\eta)}\,d\eta\;\tilde d^{\,T}((X_1,w_1),(X_2,w_2))\\
 &&=({L_2}/{T})\;\tilde d^{\,T}((X_1,w_1),(X_2,w_2)).
\end{array}
\]
A similar estimate is valid for $u<0$. Thus,
$$d^{\,T}(S_2(X_1,w_1),S_2(X_2,w_2))\le({L_2}/{T})\,\tilde d^{\,T}((X_1,w_1),(X_2,w_2)).$$
 Above arguments imply that if $T>L_2$ and $\bar k$ is small enough,
 then $S$ is a contraction in $\tilde\mathcal{B}$ w.r. to the metric defined above.
 By the Banach fixed point theorem, for any $\psi=(X_1(u),w_1(u))$ in $\tilde\mathcal{B}$
the iterated function sequence $\psi, S(\psi), S(S(\psi)), \dots$ converges uniformly on $I$
to the unique fixed point $(X(u),\;w=v_0'(u))$ of $S$ in $\tilde\mathcal{B}$.
Since $\psi=(\frac{\sqrt{\alpha^2+1}}{\alpha}\,u,0)\in\tilde\mathcal{B}$,
this point can be chosen as starting one in the iterative process.
\hfill$\square$

\subsection{Auxiliary results}
 \label{sec:axresults}

We consider a first order quasilinear system of PDE's,
$n$ equations in $n$ unknown functions $u = (u_1,\ldots, u_n)$ and two variables $x,y\in\R$,
\begin{equation}\label{E-PDE-1}
 d u/dy + A(x,y, u)\,d u/dx = b(x,y, u),
\end{equation}
where $A=(a_{ij}(x,y, u))$ is an $n\times n$ matrix, $b=(b_{i}(x,y, u))$ is an $n$-vector.

The \textit{Cauchy problem} for
(\ref{E-PDE-1})
is the problem of finding $u$ such that (\ref{E-PDE-1}) and $u(x,0)=u_0(x)$ are satisfied, where $u_0$ is given.
 When the coefficient matrix $A$ and the vector $b$ are functions of $x$ and $y$ only, the system is \textit{linear}.
When $A$ and $b$ are functions of $x$, $y$ and $u$, the system is \textit{quasilinear}.

The system (\ref{E-PDE-1}) is called \textit{hyperbolic} in the $y$-direction at $(x,y, u)$ (in an appropriate domain of the arguments of $A$) if the (right) eigenvectors of $A$ are real and span $\R^n$. In this case, let $R=[r_1,\ldots, r_n]$ be the matrix of the (right) eigenvectors $r_i$ of $A$.
For a solution $u(x,y)$ to (\ref{E-PDE-1}), the corresponding eigenvalues $\lambda_i(x,y, u)$ are called the \textit{characteristic speeds},
 the vector field $\partial_y+\lambda_i\partial_x$ is the $i$-\textit{characteristic field}, and its integral curves are $i$-\textit{characteristics}.

\vskip2mm\hskip-2mm
\textbf{Theorem~A} (see \cite{hw}\label{Th-EU1}){\it
Let the \underline{quasi-linear} system of PDE's (\ref{E-PDE-1}) be such~that

\hskip-3mm
(i) it is hyperbolic in the $y$-direction in
 $\Omega=\{|x|\le a, 0\le y\le s, \|u\|_{\infty}\le r\}$

 for some $s,r>0$;

\hskip-3mm
(ii) the matrices $A, R$ and the vector $b$ are $C^{1}$-regular in $\Omega$;

\hskip-3mm
(iii) it is satisfied an initial condition
\begin{equation}\label{E-pdes-2}
  u(x,0)= u_0(x),\quad -a\le x\le a
\end{equation}

for which $\|u_0\|_{[-a,a]}<r$ and $u_0$ is $C^{1}$-regular in $[-a,a]$.

\noindent
Then there is $\eps\in(0,s]$ such that (\ref{E-PDE-1}) and (\ref{E-pdes-2})
admit a unique $C^{1}$-regular solution $u(x,y)$ in the trapeze
 $\Pi_{K,\eps}=\{(x,y):\ |x| + K y\le a,\ 0\le y\le\eps\}$,
where $K=\max\limits\{|\lambda_i(x,y, u)|:\ (x,y, u)\in\Omega,\ 1\le i\le n\}$}.

In Theorem~A, one may use the norm
$\|u\|_{c,\infty}=\max\limits_{1\le i\le n}c_i|u_i|$ for $c_i>0$.
To~show this one should replace unknown functions $v_i= u_i/c_i$ (hence
$\|v\|_{\infty}=\|u\|_{c,\infty}$) to reduce to original Theorem~A for~$v$.
 We use Theorem~A in the proof of Theorem~\ref{T-main0} for diagonal matrices $A=R$.

\vskip1mm
The next proposition is known. For convenience of a reader we prove~it.

\begin{proposition}\label{P-twoODEs}
 Let vector functions $P$ and $Q$ satisfy the Lipschitz condition
\[
 \|P(t, u)-P(t, v)\|_{\infty}\le\bar L\|u- v\|_{\infty},\ \
 \|Q(t, u)-Q(t, v)\|_{\infty}\le\bar L\|u- v\|_{\infty},
\]
(with the same $\bar L$)
for $t\in[0,h]$ and $u, v\in\Omega\subset\R^n$ ($\Omega$ a domain).
 Let $y(t)\ (y(0)=y_0\in\Omega)$ and $z(t)\ (z(0)=z_0\in\Omega)$ are solutions to ODE's $y\,'(t) = P(t, y(t))$, and $z\,'(t) = Q(t, z(t))$, resp.,
 where $t\in[0,h]$ and $y(t),\,z(t)\in\Omega$.
 Then
 $\|y - z\|_{\,[0,h]}\le (mh + \|y_0 - z_0\|_{\infty})\,e^{\bar L\,h}$,
 where $m = \|P - Q\|_{\,[0,h]\times\Omega}$.
\end{proposition}

\textbf{Proof}. We present BVP equivalently in the integral form
\[
\begin{array}{c}
 (y\,{-}\,z)\,{-}\,(y_0\,{-}\,z_0)=\hskip-2pt\int\limits_0^t
 [P(x, y(x))\,{-}\,P(x, z(x))\,{+}\,P(x, z(x))\,{-}\,Q(x, z(x))]\,dx.
\end{array}
\]
Hence
 $\|y - z\|_{\infty}\le (\|y_0 - z_0\|_{\infty}+mh)
 +\bar L\int\nolimits_0^t\|y - z\|_{\infty}\,dx$.
From the Gronwall-Bellmann integral inequality
\[
\begin{array}{c}
 u(t)\le A +\int\limits_0^t u(x)v(x)\,dx\ (u,v>0,\ A\ge0)
 \ \Rightarrow \ u(t)\le A e^{\,\int\nolimits_0^t v(x)\,dx},
\end{array}
\]
with $A=\|y_0 - z_0\|_{\infty},\, u=\|y - z\|_{\infty}$
and $v=\bar L$ it follows the claim.\hfill$\square$

\begin{lemma}\label{appl}
Let a function $u\,{\ge}\,0$ of class $C^0([0,a])$ obeys the inequality
 $u(y)\le A\int\limits_0^y(1+u^2(\eta))^{3/2}\,d\eta$
with $A\in(0,1/a)$. Then $u(y)\le\frac{A\,y}{(1-A^2y^2)^{1/2}}$ for $y\in[0,a]$.
\end{lemma}

\textbf{Proof}. Denote $f(u)=A(1{+}u^2)^{3/2}$.
Then $U:=\frac{Ay}{(1-A^2y^2)^{1/2}}$ is the solution to the ODE $\frac{du}{dy}=f(u)$ with the initial condition $U(0)=0$. Since $A\in(0,1/a)$, $U$ is continuous in $[0,a]$.
Clearly, $U(y)$ satisfies in $[0,a)$ the integral equation
\begin{equation}\label{inteqU-3}
\begin{array}{c}
 U(y)=\int\limits_0^yf(U(\eta))\,d\eta.
\end{array}
\end{equation}
Denote by $Y(y)=U(y)-u(y)$ and $\Lambda(y)= \int_0^1\partial_u f(tU(y)+(1-t)u(y))\,dt$.
Hence $f(U(y))-f(u(y))=\Lambda(y)Y(y)$.
Since $\partial_uf(u)$ is non-negative and continuous in
$[0,\infty)$, the function $\Lambda(y)$ is also non-negative and continuous in $[0,a]$.
Furthermore, by conditions of the lemma
and (\ref{inteqU-3}), $Y(y)$ satisfies in $[0,a]$ the integral inequality $Y(y)\ge\int_0^y\Lambda(\eta)Y(\eta)\,d\eta$,
which can be written as the equation
$Y(y)=\int_0^y\Lambda(\eta)Y(\eta)\,d\eta+\phi(y)$, where $\phi(y)\in C^0([0,a])$ is
non-negative. The above linear integral
equation of Volterra type can be solved by the iterative method
$Y_0(y)=\phi(y),\ Y_{n+1}(y)=\phi(y)+\int_0^y\Lambda(\eta)Y_n(\eta)\,d\eta$.
Since all the functions of the sequence are non-negative in $[0,a]$,
their limit $Y(y)$ is also non-negative, that is
the desired estimate is valid.\hfill$\square$

\vskip2mm
The \textit{covariant derivative} of a $(0,1)$-tensor $(\mu_i)$
in $(\bar M^3,\bar g)$ is defined by
\begin{equation}\label{E-covder}
 \bar\nabla_i \,\mu_j = \mu_{\,, x_i x_j} -\sum\nolimits_k\bar\Gamma^k_{ij}\mu_{k}
\end{equation}
where
$\bar\Gamma^k_{ij}=\sum\nolimits_s\bar g^{sk}(\bar g_{is,x_j}+\bar g_{js,x_i}
-\bar g_{ij,x_s})$ are Christoffel symbols.

\begin{lemma}\label{L-EFG-1}
Let $M^2: x_3 = f(x_1,x_2)$ be the graph of a function $f\in C^2(\Pi)$
in $(\bar M^3,\bar g)$.
Then the unit normal $(n_i)$ to $M^2$, the coefficients of the 1-st and the 2-nd fundamental forms of $M^2$ are
\begin{subequations}
\begin{eqnarray}\label{E-n123}
 n_1\eq\frac{(\hat g_{12}q-\hat g_{22}p)\,\hat g_{33}}{\sqrt{EG-F^2}\,\sqrt{\det\hat g}},\ \
 n_2=\frac{(\hat g_{12}p-\hat g_{11}q)\,\hat g_{33}}{\sqrt{EG-F^2}\,\sqrt{\det\hat g}},\ \
 n_3=\frac1{\delta\,\hat g_{33}},\qquad\\
\label{E-EFG}
 E\eq \hat g_{11}+\hat g_{33}p^2,\quad
 F= \hat g_{12}+\hat g_{33}p\,q,\quad
 G= \hat g_{22}+\hat g_{33}q^2,\\
\label{E-LMN}
 L\eq (f_{x_1 x_1} {+} L_1)/\delta,\ \
 M =  (f_{x_1 x_2} {+} M_1)/\delta,\ \
 N =  (f_{x_2 x_2} {+} N_1)/\delta,\quad
\end{eqnarray}
\end{subequations}
where
$f_{x_1}=p,\ f_{x_2}=q$,
$\delta=\sqrt{(EG-F^2)/\det\hat g}\ge 1/\sqrt{\hat g_{33}}$, and
\begin{equation}\label{E-L1-M1-N1}
\hskip-2mm
\begin{array}{c}
 L_1=\hat\Gamma^3_{11}{+}2\,\hat\Gamma^3_{13}p{+}\hat\Gamma^3_{33}p^2
  {+}\delta\sum\limits_{i,j\le2}\hat g_{ij}n_j(\hat\Gamma^i_{11}{+}2\,\hat\Gamma^i_{13}p
  {+}\hat\Gamma^i_{33}p^2) {-}\Gamma^1_{11}p{-}\Gamma^2_{11}q,\\
 M_1=\hat\Gamma^3_{12}{+}\hat\Gamma^3_{23}p{+}\hat\Gamma^3_{13}q{+}\hat\Gamma^3_{33}pq
  {+}\delta\hskip-1mm\sum\limits_{i,j\le2}\hat g_{ij}n_j
  (\hat\Gamma^i_{12}+\hat\Gamma^i_{23}p+\hat\Gamma^i_{13}q
  +\hat\Gamma^i_{33}pq)\\
    -\Gamma^1_{12}p-\Gamma^2_{12}q,\\
 N_1= \hat\Gamma^3_{22}{+}2\,\hat\Gamma^3_{23}q{+}\hat\Gamma^3_{33}q^2
  {+}\delta\sum\limits_{i,j\le2}\hat g_{ij}n_j(\hat\Gamma^i_{22}{+}2\,\hat\Gamma^i_{23}q
  {+}\hat\Gamma^i_{33}q^2){-}\Gamma^1_{22}p{-}\Gamma^2_{22}q.\\
\end{array}
\end{equation}
\end{lemma}

The proof of Lemma~\ref{L-EFG-1} is based on the following

\vskip2mm
\textbf{Proposition A} (see \cite{Am2001}). {\it
The equations $x_i=\bar f_i(u_1,u_2), (i=1,2,3)$ define a regular surface $M^2$
in $(\bar M^3,\bar g)$ if and only if $\bar f_i$ are regular (of class $C^2$), and the rank
of $(\bar f_{i\,,x_j})$ is equal to 2. The first $(g_{ij})$ and the second $(b_{ij})$ fundamental forms of $M^2$ with the unit normal $(n_s)$ are given by
\begin{equation}\label{E-12Aminov}
 g_{ij}=\sum\nolimits_{\mu,\nu}\bar g_{\mu\,\nu}\bar f_{\mu\,,x_i}\bar f_{\nu\,,x_j},\quad
 b_{ij}\,n_s=\bar f_{s,\,i j}
 +\sum\nolimits_{\mu,\nu}\bar\Gamma^s_{\mu\nu}\bar f_{\mu\,,x_i}\bar f_{\nu\,,x_j},
\end{equation}
where $\bar\Gamma^s_{\mu\nu}$ are Christoffel symbols of the 2-nd kind on $\bar M^3$.}

\vskip2mm
\textbf{Proof of Lemma~\ref{L-EFG-1}}.
 From the definition
$E = \hat g(e_1,e_1), F = \hat g(e_1,e_2)$, $G = \hat g(e_2,e_2)$
it follows~(\ref{E-EFG}).
 Let $n=n_1\hat e_1{+}n_1\hat e_2{+}n_3\hat e_3$ be a unit normal to~$M^2$.
We find $n_3$ from
$
 n_3\hat g_{33}=\bar g(\hat e_3, n)=\frac{\det(e_1,e_2,\hat e_3)}{\det(e_1,e_2, n)} =\frac{\det(\hat e_1,\hat e_2,\hat e_3)}{\sqrt{EG-F^2}}=\frac1\delta.
$
Here
 $\delta^{\,2}=\frac{EG-F^2}{\det\hat g}\ge\frac1{\hat g_{33}}$,
see (\ref{E-alpha-3})$_3$.

The expressions for $n_1$ and $n_2$ of (\ref{E-n123}) follow from the linear system
$$
 \bar g(e_1, n)=n_1\hat g_{11}+n_2\hat g_{12}+n_3\hat g_{33}p=0,\ \
 \bar g(e_2, n)=n_1\hat g_{12}+n_2\hat g_{22}+n_3\hat g_{33} q=0.
$$
 From (\ref{E-12Aminov})$_2$, and $\bar g_{13}=\bar g_{23}=\bar f_{a\,,ij}=0\ (a=1,2)$, $\sum\nolimits_{a,b}\hat g_{ab}n_a n_b=1$,
we have
\begin{eqnarray*}
 b_{ij}\eq \sum\nolimits_{a,b}\hat g_{ab}n_b(b_{ij}n_a)
 =\sum\nolimits_{a,b}\hat g_{ab}n_b(\bar f_{a\,,ij}
 +\sum\nolimits_{\mu\nu}\hat\Gamma^a_{\mu\nu}\bar f_{\mu\,,i}\bar f_{\nu\,,j})\\
 \eq \hat g_{33}n_3(\bar f_{3\,,ij}
 +\sum\nolimits_{\mu\nu}\hat\Gamma^3_{\mu\nu}\bar f_{\mu\,,i}\bar f_{\nu\,,j})
 +\sum\nolimits_{a,b\le2}\hat g_{ab}n_b
 (\sum\nolimits_{\mu\nu}\hat\Gamma^a_{\mu\nu}\bar f_{\mu\,,i}\bar f_{\nu\,,j}).
\end{eqnarray*}
Hence, the coefficients $L=b_{11},M=b_{12}=b_{21}$ and $N=b_{22}$ of II are given by
\begin{eqnarray*}
\begin{array}{ccc}
 \delta\,L=&& f_{,11}+\hat\Gamma^3_{11}+2\,\hat\Gamma^3_{13}p+\hat\Gamma^3_{33}p^2
  +\delta\sum\limits_{i,j\le2}\hat g_{ij}n_j(\hat\Gamma^i_{11}+2\,\hat\Gamma^i_{13}p
  +\hat\Gamma^i_{33}pq),\\
 \delta\,N=&& f_{,22}+\hat\Gamma^3_{22}+2\,\hat\Gamma^3_{23}q+\hat\Gamma^3_{33}q^2
  +\delta\sum\limits_{i,j\le2}\hat g_{ij}n_j(\hat\Gamma^i_{22}+2\,\hat\Gamma^i_{23}q
  +\hat\Gamma^i_{33}q^2),\\
 \delta\,M=&& f_{,12}{+}\hat\Gamma^3_{12}{+}\hat\Gamma^3_{23}p{+}\hat\Gamma^3_{13}q
 {+}\bar\Gamma^3_{33}pq
 +\delta\sum\limits_{i,j\le2}\hat g_{ij}n_j(\hat\Gamma^i_{12}{+}\hat\Gamma^3_{33}pq\\
 &&{+}\hat\Gamma^i_{23}p{+}\hat\Gamma^i_{13}q{+}\hat\Gamma^i_{33}pq)
\end{array}
\end{eqnarray*}
where $f_{,\,i j}=f_{x_i x_j} -\Gamma^1_{ij}p-\Gamma^2_{ij}q$
are the covariant derivatives, see (\ref{E-covder}).
\hfill$\square$

\section{Appendix: Parallel curved surfaces}
\label{sec:appendix}

We survey basic properties of PC surfaces in aim to illustrate that the PC surfaces provide a special class of solutions to the geometrical problem.

\subsection{PC surfaces in $\R^3(k)$}
\label{sec:Ex-003}

(a) For $\alpha=\mbox{const}>0$, a solution to (\ref{E-bc0-BB2},b), (\ref{E-system-4PDEs})
 is a PC surface in $\R^3$.
 If $\bar k_2=0$, we get a cylinder $M^2: z=\sqrt{1/k^2 -(y+\alpha x)^2/(1+\alpha^2)}$ of radius $1/\bar k_1$ with the axis $\omega=(-1,\alpha,0)$.
We will build a PC surface with $c_3\ne\mbox{const}$, see Corollary~\ref{T-main3}.
Let $M_1: X^2+Z^2=R^2(Y)$ be a surface of revolution in ${\mathbb R}^3$,
where $R\ge0$ is an increasing $C^1$-regular function. Revolving about $z$-axis,
 $X=\frac{\alpha x+y}{\sqrt{1+\alpha^2}},\,
  Y=\frac{\alpha y-x}{\sqrt{1+\alpha^2}},\, Z=z,$
and replacing the function
 $\,R(\frac{t}{\sqrt{1+\alpha^2}})=\frac{r(t)}{\sqrt{1+\alpha^2}}$, we obtain
$
 z=\sqrt{R^2(Y)-X^2}={(1+\alpha^2)}^{-1}\sqrt{r^2(\alpha y-x)-(\alpha x+y)^2}.
$
A parallel $\{Y=c\}$ of above $M^2$ lies in the plane $\alpha y-x=c$, and
projects onto $xy$-plane as a line segment.

(b) Consider sphe\-rical coordinates $(\rho,\,\varphi,\,\theta)$ in the domain
$U=\{|\rho-1|\le a_1,\,|\varphi|\le a_2,\,|\theta-\pi/2|\le a_3\}$
of $\R^3$, where $0<a_1<1$ and $0<a_2<\pi$ and $0<a_3<\pi/2$.
The~\textit{curvilinear projection} onto $\Pi=\{|\rho-1|\le a_1,\,|\varphi|\le a_2,\,\theta=\frac\pi2\}$ (with $\lambda=0$) is given by
$\pi(\rho,\,\varphi,\,\theta)=(\rho,\,\varphi,\,\pi/2)$.

\vskip.5mm
(b)$_1$ Denote  $\gamma=\{|\rho-1|\le a_1,\,\varphi=0,\,\theta=\pi/2\}$ the line segment in $\Pi$. Let $\bar k_{1},\bar k_{2}$ be the functions of class $C^{1}([-a_1, a_1])$
and $l$ a vector field of class $C^{2}(\Pi)$ that is transversal but not orthogonal
to $\gamma$.
By Theorem~\ref{T-main0}, if $\bar k_{i}$ are small enough in the $C^0$-norm, then \textbf{Problem~\ref{Prob-1}} admits a unique smooth solution $M^2: \theta=f(\varphi,\rho)$
on $\Pi_{K,\eps}\,{=}\,\{|\rho-1|+K\varphi\le a_1,\,0\le\varphi\le\eps,\,\theta=\pi/2\}$.

\vskip.5mm
(b)$_2$ Assume that $M^2$ is a PC surface relative to the sphere $\beta=\{\rho=\rho_0\}$.
Take $\gamma=\{|\rho-1|\le a_1,\,\varphi=b\,(\rho-1),\,\theta=\frac\pi2\}$ for some
$b\in(0, \frac{a_2}{1-a_1})$.
The $k_1$-curvature lines of $M^2$ project onto concentric circles $\{\rho=c\}$ on $\Pi$,
hence $l=\partial_{\varphi}$ is transversal but not orthogonal to $\gamma$
(Theorem~\ref{T-main0} is applicable).
Now let $\bar k_{1}\in C^{1}([-a_1,a_1])$ and $\bar k_{2}\in C^{0}([-a_1,a_1])$ are small enough in the $C^0$-norm.
Follow the proof of Theorem~\ref{T-RPC2} (Section~\ref{sec:Th2-3}), one may show that $M^2: \theta=f(\varphi,\rho)$ can be recovered over a curvilinear rectangle $\Pi(a_1)$.

The geometric construction of a PC surface $M^2$ is as follows.
The spheres $S^2(c)=\{\rho=c\}$ intersect $M^2$ transversally by $k_1$-curvature lines.
Let $\pi_1:\R^3\setminus\{0\}\to S^2(1)$ be the \textit{radial projection} onto the unit sphere, i.e., $\pi_1(x)=x/\|x\|$.
Let $\gamma(t)$ be a $k_1$-curvature line on $M^2$, and $\gamma(t)$ belongs to $S^2(c)$ for some $c$. The curve $\gamma_0=\pi_1(\gamma)$ is homothetic  to $\gamma$
(the coefficient of homothety is $1/c$).
The great circles on $S^2(1)$ orthogonal to $\gamma_0$ and the curves
of constant distance to $\gamma_0$ form a {\em semi-geodesic net} on $S^2(1)$ near $\gamma_0$, see Lemma~\ref{L-4nets}.

\vskip1.0mm
A 1-parameter family of geodesics and their orthogonal curves on $\beta$ is called a \textit{semi-geodesic~net} (it is uniquely determined by the base curve $\gamma_0$).

\begin{figure}[ht]
\begin{center}
\includegraphics[scale=.4,angle=0,clip=true,draft=false]{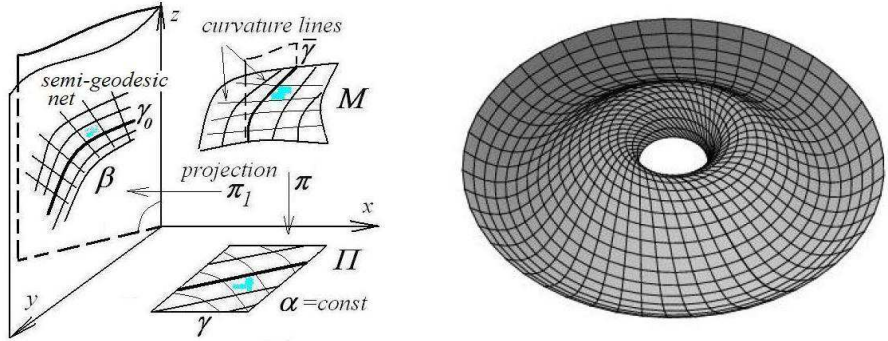}
\caption{\small (a) Projection of a PC surface $M^2$.
 \ (b) PC surface in $\R^3$ of type (c).}
\label{F-00y}
\end{center}
\end{figure}

\begin{lemma}\label{L-4nets}
There are (locally) four types of semi-geodesic nets on $(\beta, g_k)$:

\hskip-3mm
 (a) \underline{cartesian net}, $k=0,-1$: $\gamma_0$ is a line for $k=0$,
 (horocycle for $k=-1$),

\hskip-3mm
 (b) \underline{polar net}, $k=0\pm1$: $\gamma_0$ is a circle,

\hskip-3mm
 (c) \underline{evolvent net}, $k=0,\pm1$: normals to $\gamma_0$ are tangent to a curve~$\gamma_{1}$,

 ($\gamma_0$ is \textit{evolvent} of $\gamma_{1}$).

\hskip-3mm
 (d) \underline{super-parallel net}, $k=-1$: $\gamma_0$ is a line.

\noindent
 The cartesian and polar nets
 correspond to cylindrical surfaces ($k_2=0$) and  surfaces of revolution ($k_1{=}\,\mbox{const}$ along $\mathcal{F}_1$-curves,
 the axis is orthogonal to $\beta$), resp.
 (Case (c) appears on PC surfaces illustrated in Fig.~\ref{F-00y}(b)).
\end{lemma}

\textbf{Proof}. It is known that
the normals to a regular curve in $\R^2(k)$ (locally) form one of four families: (super-)parallel lines, lines through a point and enveloping a smooth curve~$\gamma_1$.\hfill$\square$

\begin{proposition}\label{L-RPC-1a}
Let $M^{2}\subset\R^3(k)$ be a PC surface-graph related to a totally umbilical surface $\beta$, and $\calf_i\ (i=1,2)$ the $k_i$-curvature lines. Then

(i) the principal curvature $k_2$ is constant along the curves of $\calf_1$,

(ii) the curves of $\calf_2$ are
geodesics, they belong to planes orthogonal to~$\beta$,

(iii) both families of curves project onto $\beta$ as a semi-geodesic net.
\end{proposition}

\textbf{Proof}. The curves of $\calf_1$ belong to totally umbilical surfaces $\beta_d$ (on the distance $d$ to $\beta$), and the curves are parallel on $M$. Hence, $\calf_2$ (that is orthogonal to $\calf_1$) consists of geodesics of $M^2$.
Let $X_i\ (i=1,2)$ be unit vector fields tangent to $\calf_i$,
$n$ a unit normal to $M^2$, $\partial_t$ a unit normal to $\beta$. Then $X_1$ is orthogonal to $n$ and $\partial_t$.
By~(2) of Lemma~\ref{L-chen1} (see Section~\ref{sec:warped2}) and Rodrigues theorem (see \cite{Top}),
$$
 X_1(\langle n, \partial_t\rangle)
 =\langle\nabla_{X_1} n, \partial_t\rangle+\langle n, \nabla_{X_1} \partial_t\rangle
 =\langle k_1 X_1, \partial_t\rangle
 +\langle n, (\log\phi)' X_1\rangle=0
$$
where $\nabla$ is the covariant derivative.
Hence, the angle between surfaces $M^2$ and $\beta_d$ along the curvature lines $\calf_1$ (the intersection) is constant.
The projections of $\calf_1$ onto $\beta$ are parallel curves $\tilde\calf_1$,
hence their orthogonal trajectories $\tilde\calf_2$ are geodesics on $\beta$.
Thus $(\tilde\calf_1,\tilde\calf_2)$ is a semi-geodesic net on $\beta$.
In~coordinates of curvature lines we have $k_{2,1}=(k_1-k_2)\frac{g_{22,\,1}}{2 g_{22}}$, see Remark~\ref{R-petcod}. Since $g_{22}=1$ ($\calf_2$-curves are unit speed geodesics), we obtain $k_{2,1}=0$, hence $k_{2}=\mbox{const}$ along $\calf_1$-curves.
One may show that (as in $\R^3$, see \cite{Ando2004}, \cite{Ando2007} and (\ref{E-k12})$_2$ in what follows)
the curves of $\calf_2$ are congruent in $\bar M^3(k)$ each to another, and lie in planes
through geodesics $\tilde\calf_2\subset\beta$ and orthogonal to~$\beta$.\hfill$\square$

\subsection{PC surfaces in a Riemannian warped product 3-space}
\label{sec:warped2}

\vskip1.0mm
Let $(S,g_k)$ be a Riemannian 2-space of constant curvature~$k$, and $\psi:I\to \R_+$.
 The \textit{Riemannian warped product} 3-space is
 $\R^3(k,\psi)=(I\times S, g^k_\psi)$ where $g^k_\psi=dt^2+\psi^2(t)g_k$.
 $\R^3(k,\psi)$ contains no open subsets of constant curvature if and only if
$(\log\psi)''+k/\psi^2\ne0$ on $I$, \cite{Chen_2008}.
 A surface $S(t)=\{t\}\times S$
 is a \textit{slice} of $\R^3(k,\psi)$.
 The {mean curvature vector} of $M^2\subset\R^3(k,\psi)$ is defined by
 $H=({\rm tr}\,h)/2$, where $h$ the {second fundamental form} of $M^2$.
 A surface $M^2$ is

 -- {totally geodesic} if $h=0$;

 -- {totally umbilical} if $h(X,Y)=g^k_\psi(X,Y)H\ \ (X,Y\in TM)$;

 -- $\mathcal{H}$-\textit{surface} if the
 vector field $\partial_t$ is tangent to $M^2$ at each point on $M^2$.

\noindent
We decompose a vector field $v$ on $\R^3(k,\psi)$ into a sum $V=\phi_V\partial_t+\tilde V$, where $\phi_V=g(V,\partial_t)$ and $\tilde V$ (a vertical component) is orthogonal to~$\partial_t$.

\begin{lemma}[\cite{Chen_2008}]\label{L-chen1}
The connection and the curvature of $\R^3(k,\psi)$ satisfy
 \begin{eqnarray*}
 (1) &&\nabla_{\partial_t}\partial_t= 0,\quad
 (2) \ \nabla_{\partial_t} X=\nabla_{X}\partial_t=(\log\psi)'X,\\
 (3) &&g(\nabla_{X}Y,\partial_t)= -g(X,Y)(\log\psi)',\
 (4) \ \nabla_{X}Y\ \mbox{is the lift of } \ \nabla^S_{X}Y  \mbox{ on } S,\\
 (5) &&\hskip-6mm R(\partial_t, X)\partial_t = (\psi''/\psi)X,\
 R(X,\partial_t)Y = \langle X, Y\rangle(\psi''/\psi)\,\partial_t,\
 R(X, Y)\partial_t = 0,\\
 &&\hskip-6mm R(X, Y) Z =(k-(\psi')^2)/\psi^2)\{\langle Y, Z\rangle X -\langle X, Z\rangle Y\}
 \quad\mbox{for}\ X,Y,Z\in TS.
\end{eqnarray*}
\end{lemma}

\begin{lemma}
 An $\mathcal{H}$-surface $\Pi_\gamma=I\times\{\gamma\}$ over a smooth curve $\gamma\subset S$ is

 (i) a ruled surface with rulings $I\times\{s\}\ (s\in\gamma)$,

 (ii) a totally geodesic in $\R^3(k,\psi)$ if and only if $\gamma$ is a geodesic in $S$.
\end{lemma}

\textbf{Proof}.
Let $h$ be the second fundamental form of $\Pi_\gamma$.
Denote $X$ the (unit) velocity field of a geodesic $\gamma$.
Using Lemma~\ref{L-chen1}, we have  on $\Pi_\gamma$:

\hskip-3.5mm
by (1): $h(\partial_t,\partial_t)=0$.
Hence $I\times\{s\}$ are rulings (geodesics in $\R^3(k,\psi)$);

\hskip-3.5mm
by (2): $\nabla_{\partial_t} X\in TM^2$, hence $h(\partial_t,X)=0$;

\hskip-3.5mm
by (4): $h(X,X)=0$ if and only if $\gamma$ is a geodesic in $S$.

\noindent
We conclude that $h=0$ when $\gamma$ is a geodesic in $S$.
On the other hand, by (3) and (4) of Lemma~\ref{L-chen1}, $\nabla^S_X X=0$ if and only if $h(X,X)=0$. Hence, if $h=0$ then $\nabla^S_X X=0$, that is $\gamma$ is a geodesic in $S$. \hfill$\square$

\vskip2mm
An $\mathcal{H}$-surface $\Pi_\gamma=I\times\{\gamma\}$ is totally umbilical with $\nabla^\perp H=0$ if and only if $\Pi_\gamma$ is totally geodesic, see \cite{Chen_2008}.
Any such $\Pi_\gamma$ over an $S$-geodesic~$\gamma$ will be named $\mathcal{H}$-\textit{plane}.
By Lemma~\ref{L-chen1}, the gaussian curvature of $\Pi_\gamma$~is $K= \psi''/\psi.$

A surface $M^{2}\subset\R^3(k,\psi)$ is called \textit{parallel curved} (PC)
relative to $S$ if it does not belong to a slice, and at each point $x\in M^2$
at least one principal direction is tangent to $S(t)$ passing through $x$. A PC surface
is \textit{regular} if such principal directions form a 1-dimensional foliation ($\calf_1$).

Proposition~\ref{L-RPC-1a} can be extended as follows

\begin{proposition}\label{L-RPC-1}
 Let $M^{2}\subset\R^3(k,\psi)$ be a regular PC surface-graph over domain in~$S$.
 Then the 2-nd family of curvature lines ($\calf_2$) consists of geodesics on $M^2$
 which lie in $\mathcal{H}$-planes.
 Two families ($\calf_1$ and $\calf_2$) of curves project onto $S$ as a semi-geodesic~net.
\end{proposition}

\noindent
Hence PC surfaces in $\R^3(k,\psi)$ represent a special class of solutions to \textbf{Problem~\ref{Prob-1}} for graphs over domains of $\mathcal{H}$-plane $\Pi_\gamma$
with $\gamma$ transversal to slices.

\end{document}